\documentclass[a4paper,reqno]{amsart}
\usepackage{color,latexsym,amssymb,amsfonts}
\usepackage{graphicx}
\usepackage{bm}
\usepackage[authoryear]{natbib}
\usepackage{tikz-cd}

\def\perm#1{[\![#1]\!]}
\def\sm#1{\left<#1\right>}

\def\sgn{\mbox{sgn}}
\def\leb{\mbox{Leb}}
\def\eqd{\,{\buildrel d \over =}\,}
\def\neqd{\,{\buildrel d \over \neq}\,}
\newcommand{\bds}{\begin{displaystyle}}
\newcommand{\eds}{\end{displaystyle}}

\newtheorem{theorem}{\quad\sc Theorem}
\newtheorem{corollary}{\quad\sc Corollary}
\newtheorem{remark}{\quad\sc Remark}
\newtheorem{lemma}{\quad\sc Lemma}
\newtheorem{proposition}{\quad\sc Proposition}
\newtheorem{example}{\quad\sc Example}

\newtheorem*{definition*}{\quad\sc Definition}

\begin{document}
\title{\bf Matching marginals and sums}

\author{Robert Griffiths}
\address{School of Mathematics, Monash University, Australia}
\email{bob.griffiths@monash.edu}

\author{Kais Hamza}
\address{School of Mathematics, Monash University, Australia}
\email{kais.hamza@monash.edu}



\date{}
\begin{abstract}
For a given set of random variables $X_1,\ldots,X_d$ we seek as large a family as possible of random variables $Y_1,\ldots,Y_d$ such that the marginal laws and the laws of the sums match: $Y_i\eqd X_i$ and $\sum_iY_i\eqd\sum_iX_i$. Under the assumption that $X_1,\ldots,X_d$ are independent and belong to any of the Meixner classes, we give a full characterisation of the random variables $Y_1,\ldots,Y_d$ and propose a practical construction by means of a finite mean square expansion. When $X_1,\ldots,X_d$ are identically distributed but not necessarily independent, using a symmetry-balancing approach we provide a universal construction with sufficient symmetry to satisfy the more stringent requirement that, for any symmetric function $g$, $g(Y)\eqd g(X)$.
\end{abstract}

\maketitle

\markboth{{\normalsize\sc Griffiths \& Hamza}}{{\normalsize\sc Matching marginals and sums}}

\section{Introduction}

The study of the behaviour of sums of random variables, be they dependent or independent, identically distributed or not, is a central problem in probability and statistics. The inverse problem of the decomposition of random variables has attracted the attention of many, including L\'evy, Cram\'er and Khintchine.

Decomposition refers to the problem of finding a $d$-dimensional random variable $(Y_1,\ldots,Y_d)$ such that $Y_1+\ldots+Y_d$ has a given distribution $\mu$. The classical decomposition problem requires that the random variables $Y_1,\ldots,Y_d$ are independent and is concerned with the existence of a solution. The celebrated L\'evy-Cram\'er Theorem states that the sum of two independent non-constant random variables is normally distributed if and only if each of the summands is normally distributed; that is if
\begin{enumerate}
\item[(i)] $X_1$ and $X_2$ are independent,
\end{enumerate}
then the following are equivalent
\begin{enumerate}
\item[(ii)] $X_1$ and $X_2$ are both normal;
\item[(iii)] $X_1+X_2$ is normal.
\end{enumerate}
An equivalent result exists for Poisson random variables. It is due to \citet{R1938}.

What can be said when $X_1$ and $X_2$ are allowed to be dependent? For example, do two standard normal random variables $X_1$ and $X_2$ such that $X_1+X_2$ is normal with mean 0 and variance 2 have to be independent or could one find a (non-Gaussian) dependent pair $(Y_1,Y_2)$ such that $Y_1$ and $Y_2$ are standard normal, and $Y_1+Y_2$ is normal with mean 0 and variance 2?

In this paper we show that the latter occurs. We do so by characterising all possible solutions as well as by developing a universal procedure for constructing concrete examples.

More generally, given a pair of random variables $(X_1,X_2)$ with joint density $f$, we seek a pair $(Y_1,Y_2)$ with joint density $\varphi$, such that $Y_1\eqd X_1$, $Y_2\eqd X_2$, $Y_1+Y_2\eqd X_1+X_2$ and $(Y_1,Y_2)\neqd(X_1,X_2)$. While one readily expects a solution to exist, the construction is surprisingly not immediate; matching the marginal distributions is easily achieved by a change of copula, and matching the distributions of the sums is trivially realised by adding and subtracting a given quantity, but enforcing both renders the problem a lot less straightforward.

\begin{example}[The Gaussian case]
\citet{S2013} (see Counterexample 10.5) suggests the following construction of a pair $(Y_1,Y_2)$ of dependent (but uncorrelated) standard normal random variables such that $Y_1+Y_2$ is normal with mean 0 and variance equal to 2:
\begin{equation}
\varphi(x_1,x_2) = \frac1{2\pi}\exp\left(-\frac12(x_1^2+x_2^2)\right)\left(1+\kappa x_1x_2(x_1^2-x_2^2)\exp\left(-\frac12(x_1^2+x_2^2)\right)\right).\label{stoyanov}
\end{equation}
Here $\kappa$ is any positive constant that ensures that $\varphi(x_1,x_2)\geq0$ or that
$$\kappa x_1x_2(x_1^2-x_2^2)\exp\left(-\frac12(x_1^2+x_2^2)\right)\geq-1.$$
This is the case for any $\kappa\leq e^2/8$, as can be easily checked.
\end{example}

While the above example provides an answer to the Gaussian case, it does not shed any light on whether other solutions exist and how to construct them, nor does it fulfil the aspiration to characterise all constructions or that to extend the problem beyond the Gaussian and independent case.

In this paper, we propose to answer, in considerable generality, the following question. Given $X=(X_1,\ldots,X_d)$, find $Y=(Y_1,\ldots,Y_d)$ such that the marginal laws of $Y$ coincide with those of $X$, and $Y_1+\ldots+Y_d\eqd X_1+\ldots+X_d$, and do so in as generic a way as possible.

After listing a few basic facts that guide our discovery, we develop in Section \ref{secmeixner} a framework for identifying all solutions when the random variables $X_1,\ldots,X_d$ are independent and have laws within a Meixner class (Theorem \ref{theorem:2}). This approach relies on a mean square expansion of the Radon-Nikodym derivative of the law of $Y$ with respect to the law of $X$, as well as on an additive property of a system of Sheffer polynomials, both of which are realised in the context of the Meixner family -- see Proposition \ref{Prop:1multi}. This ultimately leads to a recipe for constructing solutions via finite mean square expansions. It also enables us to get a factorization theorem similar to that of L\'evy-Cram\'er and Raikov when $X_1+X_2$ belongs to a Meixner class (Theorem \ref{thm:1}).

Section \ref{seccopula} looks at the dependent and identically distributed case. It uses a symmetry-balancing approach that delivers sufficient symmetry to satisfy the more stringent requirement that, for any symmetric function $g$, $g(Y)\eqd g(X)$. The construction is universal in that it applies to any $X=(X_1,\ldots,X_d)$ as long $X_1,\ldots,X_d$ admit a joint density $f$. It introduces a family of maps indexed by functions $\varepsilon$ and $\gamma$:
\begin{displaymath}
  \begin{tikzcd}
    f \arrow[r, "\mathcal{D}"] \arrow{d} & (c,\phi) \arrow{d}{T_{\varepsilon,\gamma}\otimes\mathrm{id}} \\
    \varphi \arrow[r, leftarrow, "\mathcal{R}"] & (\theta_{\varepsilon,\gamma},\phi) \\
  \end{tikzcd}
\end{displaymath}
where $\mathcal{D}(f)$ is the ``disintegration'' of $f$ into its copula density $c$ and marginal law $\phi$, $\theta_{\varepsilon,\gamma}$ is a ``perturbation'' of $c$ of the form $c-\sum_\mathfrak{q}\varepsilon(\mathfrak{q})\gamma_\mathfrak{q}$, and $\mathcal{R}(\theta_{\varepsilon,\gamma},\phi)$ is the ``recombination'' of the copula density $\theta_{\varepsilon,\gamma}$ and the marginal law $\phi$ into the density $\varphi$. While the construction delivers much more than matching marginals and sums, it does so by imposing significant symmetry.

Our aim is to construct a pair $(Y_1,Y_2)$ that matches a given pair $(X_1,X_2)$ in regards to the marginal laws, $Y_1\eqd X_1$ and $Y_2\eqd X_2$, and the law of the sum, $Y_1+Y_2\eqd X_1+X_2$. We start by assuming that $X_1$ and $X_2$ are independent so that $(X_1,X_2)$ has a product law with respect to which the law of $(Y_1,Y_2)$ must be absolutely continuous. We call $H$ the Radon-Nikodym derivative.

We propose to use the method of expansion by orthogonal functions to identify such functions. Suppose the laws of $X_1$ and $X_2$ allow for complete sets of orthogonal functions $(\mathfrak{f}_{1,n})_n$ and $(\mathfrak{f}_{2,n})_n$:
$$\forall x_1,x_2,\ \mathfrak{f}_{1,0}(x_1)=\mathfrak{f}_{2,0}(x_2)=1\mbox{ and }\forall m\neq n,\ \mathbb{E}[\mathfrak{f}_{1,m}(X_1)\mathfrak{f}_{1,n}(X_1)]=\mathbb{E}[\mathfrak{f}_{2,m}(X_2)\mathfrak{f}_{2,n}(X_2)]=0,$$
so that $H$ admits a mean square expansion
$$H(x_1,x_2)=1+\sum_{n_1+n_2\geq1}H_{n_1,n_2}\mathfrak{f}_{1,n_1}(x_1)\mathfrak{f}_{2,n_2}(x_2),$$
provided $H(X_1,X_2)$ is square integrable or that
$$\sum_{n_1,n_2}H_{n_1,n_2}^2\mathbb{E}[\mathfrak{f}_{1,n_1}(X_1)^2]\mathbb{E}[\mathfrak{f}_{2,n_2}(X_2)^2]<+\infty.$$
The problem therefore reduces to finding a measurable non-negative $H$ the coefficients of which satisfy
\begin{equation}\label{GenCondExpansion}
\forall n_1,n_2,\ H_{n_1,0}=H_{0,n_2}=0,\mbox{ and }\sum_{n_1+n_2\geq1}H_{n_1,n_2}\mathbb{E}[\mathfrak{f}_{1,n_1}(X_1)\mathfrak{f}_{2,n_2}(X_2)|X_1+X_2]=0.
\end{equation}
The difficulty then resides in identifying functions $H$ for which the above holds.

In this paper, we take two distinct approaches. The first relies on an additivity property enjoyed by the so-called Meixner family that allows for a convenient evaluation of $\mathbb{E}[\mathfrak{f}_{1,n_1}(X_1)\mathfrak{f}_{2,n_2}(X_2)|X_1+X_2]$, when $X_1$ and $X_2$ are independent but not necessarily identically distributed. The other draws on added symmetry built into the problem. More specifically we shall assume that $X_1$ and $X_2$ are exchangeable, but not necessarily independent, and seek a function $H$ of the form $H(x_1,x_2)=1+h(x_1,x_2)-h(x_2,x_1)$ such that for any $x$, $\mathbb{E}[h(x,X_1)]=\mathbb{E}[h(X_1,x)]=0$. This is easily seen to be a necessary and sufficient condition to achieve the more stringent requirement that $g(Y_1,Y_2)\eqd g(X_1,X_2)$ for any symmetric (and measurable) $g$.

The problem therefore reduces to finding a non-symmetric $h$ such that
$$1+h(x_1,x_2)-h(x_2,x_1)\geq0\mbox{ and }\forall x,\ \mathbb{E}[h(X_1,x)] = \mathbb{E}[h(x,X_1)] = 0.$$

We stress that these facts alone are not sufficient to give an insight into how $h$ or $H$ can be constructed.

We end this section by observing that in the Stoyanov example \eqref{stoyanov}, $\varphi$ is precisely of this form: $\varphi(x_1,x_2) = \phi(x_1)\phi(x_2)\big(1+h(x_1,x_2)-h(x_2,x_1)\big)$.

\section{The independent case -- An expansion approach}\label{secmeixner}
\subsection{The Meixner class of distributions}
Given a law $\mu$ with finite moments of all order, a simple application of the Gram-Schmidt method enables the construction of a unique sequence of polynomials $(P_n)_{n\geq0}$ with the following properties:
\begin{enumerate}
\item the leading term of the polynomial $P_n$ is the monomial of degree $n$, $\pi_n(x)=x^n$;
\item for $m\neq n$, $$\left<P_m,P_n\right>_\mu=\int P_m(x)P_n(x)\mu(dx)=\mathbb{E}[P_m(X)P_n(X)]=0,$$
where $X$ has law $\mu$.
\end{enumerate}
The sequence of orthogonal polynomials $(P_n)_n$ must satisfy a three-term recurrence relation
\begin{equation}\label{3termrec}
xP_n(x) = A_nP_{n+1}(x)+B_nP_n(x)+C_nP_{n-1}(x),
\end{equation}
where, because of (1), we have in fact $A_n=1$, and where $B_n$ and $C_n$ are real numbers with $C_n\geq0$.

\citet{M1934} characterized those distributions for which $\sum_{n=0}^\infty P_n(x)z^n/n!$ can be written as $e^{x\mathfrak{u}(z)}/M(\mathfrak{u}(z))$, where $M$ is the moment generating function ($M(t)=\int e^{tx}\mu(dx)$), supposed to be finite in an open interval containing 0, and for a functions $\mathfrak{u}$ such that $\mathfrak{u}(z)$ has a power series expansion in $z$ with $\mathfrak{u}(0)=0$ and $\mathfrak{u}^\prime(0)=1$.

In this case,
\begin{equation}
\mathfrak{G}(z,x) = \frac{e^{x\mathfrak{u}(z)}}{M(\mathfrak{u}(z))} = \sum_{n=0}^\infty P_n(x)z^n/n!,
\label{genfn:0}
\end{equation}
is called the generating function of the law $\mu$ (or of the orthogonal polynomials $(P_n)_n$). It satisfies the property that, if $X$ has law $\mu$, then
$$\mathbb{E}[\mathfrak{G}(z,X)] = \sum_{n=0}^\infty\mathbb{E}[P_n(X)]\frac{z^n}{n!} = \sum_{n=0}^\infty\mathbb{E}[P_n(X)P_0(X)]\frac{z^n}{n!} = 1.$$

Meixner used the property of orthogonality to characterize $\mathfrak{u}$ and therefore distributions on which the polynomials are orthogonal. \citet{E1964} shows that these polynomials form a complete orthogonal system in $L^2(\mu)$.

Generating functions of the form (\ref{genfn:0}) when the polynomials are not necessarily orthogonal generate Sheffer polynomials, with the only orthogonal polynomials in the class being the Meixner polynomials. Sheffer polynomials are important in constructing martingales in L{\'e}vy processes \citep{S2000}.

Let $\mathfrak{v}$ denote the inverse of $\mathfrak{u}$. \citet{M1934} shows that $\mathfrak{v}$ necessarily solves the Ricatti differential equation with constant coefficients:
$$\mathfrak{v}'(t)=(1-\mathfrak{av}(t))(1-\mathfrak{bv}(t)),$$
with $\mathfrak{a}+\mathfrak{b}\in\mathbb{R}$ and $\mathfrak{a}\mathfrak{b}\in\mathbb{R}$. Following \cite{E1964}, we distinguish five types of distributions depending on the values of $\mathfrak{a}$ and $\mathfrak{b}$.
Each is fixed to have arbitrary mean $\mathfrak{m}$ and variance $\mathfrak{s}^2$. Below, we give expressions for $\mathfrak{u}$, $M$ and $\mathfrak{G}$ in each case. We also identify the corresponding laws and polynomials. $\mathfrak{h}_n=\mathbb{E}\big[P_n(X)^2\big]$ plays a significant role in the expansions we rely on. We give expressions for those as well.

\begin{enumerate}

\item[(I)] \underline{$\mathfrak{a}=\mathfrak{b}=0$}.
$\displaystyle\mathfrak{u}(z)=z$, $\displaystyle M(t)=e^{\mathfrak{m}t+\mathfrak{s}^2t^2/2}$, $\displaystyle \mathfrak{G}(z,x)=e^{(x-\mathfrak{m})z-\frac12\mathfrak{s}^2z^2}$, $\displaystyle\mathfrak{h}_n=\mathfrak{s}^{2n}n!$,
$\mu$ is a normal law and $(P_n)_n$ are the Hermite polynomials.

\item[(II)] \underline{$\mathfrak{a}=\mathfrak{b}\neq0$}.
$\displaystyle\mathfrak{u}(z)=z/(1-\mathfrak{a}z)$, $\displaystyle M(t)=e^{\left(\frac{\mathfrak{s}^2}{\mathfrak{a}}+\mathfrak{m}\right)t}(1+\mathfrak{a}t)^{-\frac{\mathfrak{s}^2}{\mathfrak{a}^2}}$, $\displaystyle \mathfrak{G}(z,x)=\left(1-\mathfrak{a}z\right)^{-\frac{\mathfrak{s}^2}{\mathfrak{a}^2}}e^{\left(x-\mathfrak{m}-\frac{\mathfrak{s}^2}{\mathfrak{a}}\right)\frac{z}{1-\mathfrak{a}z}}$, $\displaystyle \mathfrak{h}_n=\mathfrak{a}^{2n}n!\prod_{i=1}^n(\mathfrak{s}^2/\mathfrak{a}^2+i-1)$,
$\mu$ is a generalized gamma law and $(P_n)_n$ are the Laguerre polynomials.

\item[(III)] \underline{$\mathfrak{a}\neq\mathfrak{b}$, $\mathfrak{ab}=0$}. Suppose (wlog) $\mathfrak{b}=0$.
$\displaystyle\mathfrak{u}(z)=\big(1/\mathfrak{a}\big)\ln\big(1/(1-\mathfrak{a}z)\big)$, $\displaystyle M(t)=e^{(\mathfrak{m}+\mathfrak{s}^2/\mathfrak{a})t+\frac{\mathfrak{s}^2}{\mathfrak{a}^2}\left(e^{-\mathfrak{a}t}-1\right)}$, $\displaystyle \mathfrak{G}(z,x)=e^{\frac{\mathfrak{s}^2}{\mathfrak{a}}z}(1-\mathfrak{a}z)^{\frac1{\mathfrak{a}}(\mathfrak{m}-x)+\frac{\mathfrak{s}^2}{\mathfrak{a}^2}}$, $\displaystyle \mathfrak{h}_n=\mathfrak{s}^{2n}n!$,
$\mu$ is a generalized Poisson law and $(P_n)_n$ are the Charlier polynomials.

\item[(IV)] \underline{$\mathfrak{a}\neq\mathfrak{b}$, $\mathfrak{ab}\neq0$, $\mathfrak{a},\mathfrak{b}\in\mathbb{R}$}. We distinguish two cases, $\mathfrak{a}\mathfrak{b}>0$ and $\mathfrak{a}\mathfrak{b}<0$.
$\displaystyle\mathfrak{u}(z)=\big(1/(\mathfrak{b}-\mathfrak{a})\big)\ln\big((1-\mathfrak{a}z)/(1-\mathfrak{b}z)\big)$, $\displaystyle M(t)=e^{(\mathfrak{m}+\mathfrak{s}^2/\mathfrak{a})t}\left(\frac{\mathfrak{a}-\mathfrak{b}}{\mathfrak{a}-\mathfrak{b}e^{(\mathfrak{b}-\mathfrak{a})t}}\right)^{\frac{\mathfrak{s}^2}{\mathfrak{a}\mathfrak{b}}}$, $\displaystyle \mathfrak{G}(z,x)=\big(1-\mathfrak{b}z\big)^{-\frac{\mathfrak{s}^2}{\mathfrak{a}\mathfrak{b}}}
\left(\frac{1-\mathfrak{b}z}{1-\mathfrak{a}z}\right)^{\frac{\mathfrak{m}-x}{\mathfrak{b}-\mathfrak{a}}+\frac{\mathfrak{s}^2}{\mathfrak{a}(\mathfrak{b}-\mathfrak{a})}}$, $\displaystyle \mathfrak{h}_n=(\mathfrak{a}\mathfrak{b})^nn!\prod_{i=1}^n(\mathfrak{s}^2/(\mathfrak{a}\mathfrak{b})+i-1)$,
$\mu$ is a generalized negative binomial law ($\mathfrak{a}\mathfrak{b}>0$) or a generalized binomial law ($\mathfrak{a}\mathfrak{b}<0$ and $\mathfrak{s}^2/(\mathfrak{a}\mathfrak{b})$ integer), and $(P_n)_n$ are the the Meixner ($\mathfrak{a}\mathfrak{b}>0$) or Krawtchouk ($\mathfrak{a}\mathfrak{b}<0$ and $\mathfrak{s}^2/(\mathfrak{a}\mathfrak{b})$ integer) polynomials.

\item[(V)] \underline{$\mathfrak{a}\neq\mathfrak{b}$, $\mathfrak{b}=\overline{\mathfrak{a}}$}. $\displaystyle\mathfrak{u}(z)=\frac1{2\Im(\mathfrak{a})}\arg((1-\Re(\mathfrak{a})z)^2+\Im(\mathfrak{a})^2z^2-2i\Im(\mathfrak{a})z(1-\Re(\mathfrak{a})z)$ and $\mu$ is a generalised hypergeometric distribution -- see \citet{E1964} for details.

\end{enumerate}

A useful summary of the orthogonal polynomials with standard notion is in Appendix B of \citet{S2000}.

\subsection{Additivity of bivariate random variables in a Meixner class}
Independent random variables within the respective Meixner classes are additive. Indeed, let $\mathcal{M}(\mathfrak{u})$ denote the Meixner class associated with $\mathfrak{u}$. If $\mu_1,\mu_2\in\mathcal{M}(\mathfrak{u})$ and have respective moment generating functions $M_1$ and $M_2$, respective means $\mathfrak{m}_1$ and $\mathfrak{m}_2$, and respective variances $\mathfrak{s}_1^2$ and $\mathfrak{s}_2^2$, then $M=M_1M_2$ solves
$$\frac{M'(t)}{M(t)} = \frac{M_1'(t)}{M_1(t)}+\frac{M_2'(t)}{M_2(t)} = (\mathfrak{s}_1^2\mathfrak{v}(t)+\mathfrak{m}_1)+(\mathfrak{s}_2^2\mathfrak{v}(t)+\mathfrak{m}_2) = (\mathfrak{s}_1^2+\mathfrak{s}_2^2)\mathfrak{v}(t)+(\mathfrak{m}_1+\mathfrak{m}_2);$$
that is $\mu=\mu_1\ast\mu_2\in\mathcal{M}(\mathfrak{u})$.

An inspection of the five types listed above, shows that, except for the case $\mathfrak{ab}<0$, all laws within the Meixner classes are scalable. That is $\mathfrak{s}^2$ is a completely free parameter (and so is $\mathfrak{m}$). It follows that any $\mu\in\mathcal{M}(\mathfrak{u})$ can be decomposed in any number of independent and identically distributed random variables all of which belong to the same Meixner class $\mathcal{M}(\mathfrak{u})$.
When the decomposition is into $n$ random variables their law is obtained by solving
$$\frac{M'(t)}{M(t)} = \frac{\mathfrak{s}^2}{n}\mathfrak{v}(t)+\frac{\mathfrak{m}}{n}.$$
In other words, expect for the case $\mathfrak{ab}<0$, laws within the Meixner classes are infinitely divisible.

In the sequel we mostly consider random variables that belong to the same Meixner class but are not necessarily equally distributed. Such random variables differ only through their means and variances. We shall therefore qualify all objects with the parameter pair $\mathfrak{r}=(\mathfrak{m},\mathfrak{s}^2)$. We shall for example write $M(t;\mathfrak{r})$, $\mathfrak{G}(z,x;\mathfrak{r})$, $P_n(x;\mathfrak{r})$, $\mathfrak{h}_n(\mathfrak{r})$ etc whenever we consider the law within $\mathcal{M}(\mathfrak{u})$ that has a mean and a variance given by $\mathfrak{r}$. Such a law will be denoted by $\mu(\mathfrak{r})$.

A crucial step in proving the main result of this section is an expression for $\mathbb{E}\big [P_{n_1}(X_1;\mathfrak{r}_1)P_{n_2}(X_2;\mathfrak{r}_2)\mid X_1+X_2\big ]$ that allows for a simplification of \eqref{GenCondExpansion}. This is done in the next proposition and leads to a secondary result of the Cram\'er-Raikov type not directly related to the questions raised above but of intrinsic importance.

It uses the completeness of the polynomials $(P_n)_n$ to show that projections on $P_n(X)$ define the expectation conditional on $X$.

\begin{lemma}\label{Lemma:0}
Suppose $X$ belongs to one of the Meixner classes characterised by the polynomials $(P_n)_n$. Let $Y$ be square integrable. If there exists $g$ measurable such that $g(X)$ is square integrable and
$\forall n\geq0$, $\mathbb{E}[YP_n(X)] = \mathbb{E}[g(X)P_n(X)]$, then $\mathbb{E}[Y|X]=g(X)$ (a.s.).
\end{lemma}

\begin{proposition}\label{Prop:1}
Let $X_1,X_2$ be independent random variables belonging to the same Meixner class and having parameters $\mathfrak{r}_1$ and $\mathfrak{r}_2$. For $n\geq0$,
\begin{equation}
P_n(x_1+x_2;\mathfrak{r}_1+\mathfrak{r}_2)
= \sum_{k=0}^n{n\choose k}P_k(x_1;\mathfrak{r}_1)P_{n-k}(x_2;\mathfrak{r}_2)
\label{Runge:1}
\end{equation}
and, for $n_1,n_2\geq0$,
\begin{equation}
\mathbb{E}\big [P_{n_1}(X_1;\mathfrak{r}_1)P_{n_2}(X_2;\mathfrak{r}_2)\mid X_1+X_2\big ]
= {n_1+n_2\choose n_1}\frac{\mathfrak{h}_{n_1}(\mathfrak{r}_1)\mathfrak{h}_{n_2}(\mathfrak{r}_2)}{\mathfrak{h}_{n_1+n_2}(\mathfrak{r}_1+\mathfrak{r}_2)}P_{n_1+n_2}(X_1+X_2;\mathfrak{r}_1+\mathfrak{r}_2).
\label{cond:100}
\end{equation}
\end{proposition}
\begin{proof}
From the fact that, for two independent random variables that belong to the same Meixner class, say $\mathcal{M}(\mathfrak{u})$, $M(t;\mathfrak{r}_1)M(t;\mathfrak{r}_2)=M(t;\mathfrak{r}_1+\mathfrak{r}_2)$, we deduce that
$$\mathfrak{G}(z,x_1,\mathfrak{r}_1)\mathfrak{G}(z,x_2,\mathfrak{r}_2) = \frac{e^{x_1\mathfrak{u}(z)}}{M(\mathfrak{u}(z);\mathfrak{r}_1)}\frac{e^{x_2\mathfrak{u}(z)}}{M(\mathfrak{u}(z);\mathfrak{r}_2)} =
\frac{e^{(x_1+x_2)\mathfrak{u}(z)}}{M(\mathfrak{u}(z);\mathfrak{r}_1+\mathfrak{r}_2)} = \mathfrak{G}(z,x_1+x_2,\mathfrak{r}_1+\mathfrak{r}_2).$$
(\ref{Runge:1}) follows from an expansion of both sides of this identity and from identifying the coefficients of $z^n/n!$.

Using the independence of $X_1$ and $X_2$ and the orthogonality of $P_k$ and $P_\ell$, $k\neq\ell$, we deduce that
\begin{equation}
\mathbb{E}\big [P_{n_1+n_2}(X_1+X_2;\mathfrak{r}_1+\mathfrak{r}_2)P_{n_1}(X_1;\mathfrak{r}_1)P_{n_2}(X_2;\mathfrak{r}_2)\big ] = {n_1+n_2\choose n_1}\mathfrak{h}_{n_1}(\mathfrak{r}_1)\mathfrak{h}_{n_2}(\mathfrak{r}_2)
\label{exp:100}
\end{equation}
and that, for $n\neq n_1+n_2$,
\begin{equation}
\mathbb{E}\big [P_n(X_1+X_2;\mathfrak{r}_1+\mathfrak{r}_2)P_{n_1}(X_1;\mathfrak{r}_1)P_{n_2}(X_2;\mathfrak{r}_2)\big ] = 0.
\label{exp:101}
\end{equation}
Combining (\ref{exp:100}) and (\ref{exp:101}) and using Lemma \ref{Lemma:0}, we obtain \eqref{cond:100}.
%\begin{equation*}
%\mathbb{E}\big [P_{m}(X_1;\theta_1)P_{n}(X_2;\theta_2)\mid X_1+X_2\big ] = {m+n\choose n}\frac{h_m(\theta_1)h_n(\theta_2)}{h_{m+n}(\theta_1+\theta_2)}P_{m+n}(X_1+X_2;\theta_1+\theta_2).
%\end{equation*}
\end{proof}

\citet{E1964} used the Runge type identity (\ref{Runge:1}) \citep{R1914} that these polynomials satisfy
to study bivariate expansions of distributions within a Meixner class which have random elements in common.

The next theorem is a Cram\'er-Raikov-type result. It is not needed for the remainder of the paper but demonstrates the reach of the techniques used herein.
It shows that absolutely continuous and identically distributed factorizations of Meixner laws are limited to Meixner laws of the same class.

\begin{theorem}\label{thm:1}
Let $\mu(2\mathfrak{r})\in\mathcal{M}(\mathfrak{u})$. Suppose $\mu(2\mathfrak{r})=\nu\ast\nu$ with $\nu\ll\mu(\mathfrak{r})$. Suppose further that the Radon-Nikodym derivative, $d\nu/d\mu(\mathfrak{r})$, is square integrable. Then $\nu=\mu(\mathfrak{r})\in\mathcal{M}(\mathfrak{u})$.
\end{theorem}
\begin{proof}
Let $H=d\nu/d\mu(\mathfrak{r})$. Then $H$ admits a mean square expansion $\displaystyle\sum_nH_nP_n(x;\mathfrak{r})$. Suppose $(X_1,X_2)$ has law $\mu(\mathfrak{r})\otimes\mu(\mathfrak{r})$, and $X=X_1+X_2$. Note that since $\mathbb{E}[H(X_1)]=1$, $H_0=1$. It follows that $X$ has law $\mu(2\mathfrak{r})$ and if $(Z_1,Z_2)$ has law $\nu\otimes\nu$ and $Z_1+Z_2=X$, then $d(\nu\otimes\nu)/d\big(\mu(\mathfrak{r})\otimes\mu(\mathfrak{r})\big)(x_1,x_2)=H(x_1)H(x_2)$ and, for any bounded and measurable $G$,
$$\mathbb{E}[G(X)] = \mathbb{E}[G(Z_1+Z_2)] = \mathbb{E}[G(X_1+X_2)H(X_1)H(X_2)] = \mathbb{E}[G(X)H(X_1)H(X_2)].$$
Also, since $H$ admits a mean square expansion (with respect to $\mu(\mathfrak{r})$), then
$$\lim_m\mathbb{E}\Big[\Big|\sum_{n_1,n_2=0}^mH_{n_1}H_{n_2}P_{n_1}(X_1;\mathfrak{r})P_{n_2}(X_2;\mathfrak{r})-H(X_1)H(X_2)\Big|\Big]=0$$
and
$$\mathbb{E}[G(X)H(X_1)H(X_2)] = \sum_{n_1,n_2}H_{n_1}H_{n_2}\mathbb{E}[G(X)P_{n_1}(X_1;\mathfrak{r})P_{n_2}(X_2;\mathfrak{r})].$$
Note that $G$ is bounded. Now using Proposition \ref{Prop:1},
\begin{eqnarray*}
\mathbb{E}[G(X)H(X_1)H(X_2)] & = & \sum_{n_1,n_2}H_{n_1}H_{n_2}\mathbb{E}[G(X)P_{n_1}(X_1;\mathfrak{r})P_{n_2}(X_2;\mathfrak{r})]\\
& = & \sum_{n_1,n_2}H_{n_1}H_{n_2}{n_1+n_2\choose n_1}\frac{\mathfrak{h}_{n_1}(\mathfrak{r})\mathfrak{h}_{n_2}(\mathfrak{r})}{\mathfrak{h}_{n_1+n_2}(2\mathfrak{r})}\mathbb{E}[G(X)P_{n_1+n_2}(X;2\mathfrak{r})]\\
& = & \sum_n\Big(\sum_{k=0}^nH_{k}H_{n-k}\frac{\mathfrak{h}_k(\mathfrak{r})\mathfrak{h}_{n-k}(\mathfrak{r})}{k!(n-k)!}\Big)\frac{n!}{\mathfrak{h}_n(2\mathfrak{r})}\mathbb{E}[G(X)P_n(X;2\mathfrak{r})].
\end{eqnarray*}
Writing $\kappa_n$ for the coefficient of $\mathbb{E}[G(X)P_n(X;2\mathfrak{r})]$ in the above sum, we get that
$$\mathbb{E}[G(X)] = \mathbb{E}[G(X)H(X_1)H(X_2)] = \mathbb{E}\Big[G(X)\sum_n\kappa_nP_n(X,2\mathfrak{r})\Big].$$
As this must be true for any bounded measurable $G$, we deduce that $\sum_n\kappa_nP_n(X,2\mathfrak{r})=1$ or that $\kappa_n=\delta_{0n}$. Recall that $H_0=1$. As $\kappa_n$ is essentially the convolution of $H_n\mathfrak{h}_n(\mathfrak{r})/n!$ with itself, we deduce that $H_n\mathfrak{h}_n(\mathfrak{r})/n!=\delta_{0n}$ and consequently that $H_n=\delta_{0n}$. In other words, $H=1$ and $\nu=\mu(\mathfrak{r})$.
\end{proof}

\begin{example}
A direct application of Theorem \ref{thm:1} is a reduced version of the famous Cram\'er Theorem and an extension to (generalized) gamma distributions. Let $Z_1$ and $Z_2$ be independent and identically distributed random variables, and suppose they admit a (common) density.
\begin{enumerate}
\item If $Z_1+Z_2$ is normal, then so are $Z_1$ and $Z_2$.
\item If $Z_1+Z_2$ has a (generalized) gamma distribution, then so do $Z_1$ and $Z_2$.
\end{enumerate}
Similar statements can be made for the other Meixner types.
\end{example}

\subsection{Matching within the Meixner classes}

There are also non-independent random variables $Y_1,Y_2$ in the classes  such that $Y_1+Y_2$ has the same distribution as in the independent case. We characterize these distributions in an arbitrary multidimensional setting.

For $\mathfrak{r}=(\mathfrak{r}_1,\ldots,\mathfrak{r}_d)$, let $Y_1,\ldots,Y_d$ have marginal laws $\mu(\mathfrak{r}_1),\ldots,\mu(\mathfrak{r}_d)$. Their joint distribution $\mu_\bullet(\mathfrak{r})$ is absolutely continuous with respect to the product measure $\mu_{\otimes}(\mathfrak{r})=\mu(\mathfrak{r}_1)\otimes\ldots\otimes\mu(\mathfrak{r}_d)$ and the Radon-Nikodym derivative $H$ has an expansion in mean square
\begin{equation}
H(x) = \frac{d\mu_\bullet(\mathfrak{r})}{d\mu_\otimes(\mathfrak{r})}(x) = \sum_{n_1,\ldots,n_d\geq0}H_{n_1,\ldots,n_d}\prod_{i=1}^dP_{n_i}(x_i;\mathfrak{r}_i),
\label{expansion:0}
\end{equation}
provided that
\begin{equation}
\int H(x)^2d\mu_\otimes(\mathfrak{r})(x) = \sum_{n_1,\ldots,n_d\geq0}H_{n_1,\ldots,n_d}^2\prod_{i=1}^dh_{n_i}(\mathfrak{r}_i) < \infty.
\label{expansion:0-assumption}
\end{equation}
If $X=(X_1,\ldots,X_d)$ has law $\mu_\otimes(\mathfrak{r})$, then
$$\mathbb{E}\big [\prod_{i=1}^dP_{n_i}(Y_i;\mathfrak{r}_i)\big] = \mathbb{E}\big [H(X)\prod_{i=1}^dP_{n_i}(X_i;\mathfrak{r}_i)\big] = H_{n_1,\ldots,n_d}\prod_{i=1}^dh_{n_i}(\mathfrak{r}_i).$$
See \citet{L1958,L1963}. Note that $H_{0,\ldots,0}=1$ and, $Y_1,\ldots,Y_d$ are independent if and only if $H_{n_1,\ldots,n_d}=0$ whenever $(n_1,\ldots,n_d)\neq(0,\ldots,0)$. For $n=(n_1,\ldots,n_d)$, we shall write $\sm{n}$ for $n_1+\ldots+n_d$. Similarly for $x=(x_1,\ldots,x_d)$ and $\mathfrak{r}=(\mathfrak{r}_1,\ldots,\mathfrak{r}_d)$. We shall also write $H_n$ for $H_{n_1,\ldots,n_d}$ and ${\sm{n}\choose n}$ for the multinomial coefficient $(\sm{n}!)/(n_1!\ldots n_d!)$.
The next proposition is a straightforward extension of Proposition \ref{Prop:1} to a multidimensional setting.

\begin{proposition}\label{Prop:1multi}
Let $X_1,\ldots,X_d$ be independent random variables belonging to the same Meixner class and having parameters $\mathfrak{r}_1,\ldots,\mathfrak{r}_d$. For $m\geq0$,
\begin{equation}
P_m(\sm{x};\sm{\mathfrak{r}}) = \sum_{\sm{n}=m}{m\choose n}\prod_{i=1}^dP_{n_i}(x_i;\mathfrak{r}_i)
\label{Runge:1multi}
\end{equation}
and
\begin{equation}
\mathbb{E}\big [\prod_{i=1}^dP_{n_i}(X_i;\mathfrak{r}_i)\mid\sm{X}\big ]
= {\sm{n}\choose n}\frac{\prod_{i=1}^d\mathfrak{h}_{n_i}(\mathfrak{r}_i)}{\mathfrak{h}_{\sm{n}}(\sm{\mathfrak{r}})}P_{\sm{n}}(\sm{X};\sm{\mathfrak{r}}).
\label{cond:100multi}
\end{equation}
\end{proposition}

\begin{theorem}\label{theorem:2}
Suppose $Y_1,\ldots,Y_d$ belong to one of the Meixner classes and suppose they have a multivariate distribution $d\mu_\bullet(\mathfrak{r})=Hd\mu_\otimes(\mathfrak{r})$ where $H$ is given by \eqref{expansion:0} and satisfies \eqref{expansion:0-assumption}. Then
\begin{enumerate}
\item $Y_1,\ldots,Y_d$ have laws $\mu(\mathfrak{r}_1),\ldots,\mu(\mathfrak{r}_d)$ if and only if for every $i=1,\ldots,d$ and every $m\geq1$, $H_{n(i,m)}=0$, where $n(i,m)=(0,\ldots0,m,0,\ldots,0)$ and $m$ appears in position $i$.
\item $Y_1+\ldots+Y_d$ has law $\mu(\mathfrak{r}_1)\ast\ldots\ast\mu(\mathfrak{r}_d)$ if and only if for every $m\geq1$,
\begin{equation}
\sum_{\sm{n}=m}H_n\frac1{\prod_{i=1}^dn_i!}\prod_{i=1}^d\mathfrak{h}_{n_i}(\mathfrak{r}_i)=0.
\label{condition:0}
\end{equation}
\end{enumerate}
\end{theorem}
\begin{proof}

\begin{enumerate}
\item Let $X$ have law $\mu_\otimes(\mathfrak{r})$ and $Y$ have law $\mu_\bullet(\mathfrak{r})$. For $G$ bounded and measurable,
$$\mathbb{E}[G(Y_j)] = \mathbb{E}[G(X_j)H(X)] = \sum_nH_n\mathbb{E}[G(X_j)\prod_{i=1}^dP_{n_i}(X_i;\mathfrak{r}_i)]
= \mathbb{E}\Big[G(X_j)\sum_{n_j}H_{n(j,n_j)}P_{n_j}(X_j;\mathfrak{r}_j)\Big].$$
Therefore, $Y_j\eqd X_j$ if and only if for any $n_j\geq1$, $H_{n(j,n_j)}=0$.

\item We use a similar approach to that used in the proof of Theorem \ref{thm:1}. For $G$ bounded and measurable,
\begin{eqnarray*}
\mathbb{E}[G(\sm{Y})] & = & \mathbb{E}[G(\sm{X})H(X)]\\
& = & \sum_nH_n\mathbb{E}[G(\sm{X})\prod_{i=1}^dP_{n_i}(X_i;\mathfrak{r}_i)]\\
& = & \sum_nH_n{\sm{n}\choose n}\frac{\prod_{i=1}^d\mathfrak{h}_{n_i}(\mathfrak{r}_i)}{\mathfrak{h}_{\sm{n}}(\sm{\mathfrak{r}})}\mathbb{E}[G(\sm{X})P_{\sm{n}}(\sm{X};\sm{\mathfrak{r}})]\\
& = & \sum_m\left(\sum_{\sm{n}=m}H_n\frac{\prod_{i=1}^d\mathfrak{h}_{n_i}(\mathfrak{r}_i)}{\prod_{i=1}^dn_i!}\right)\frac{m!}{\mathfrak{h}_{m}(\sm{\mathfrak{r}})}\mathbb{E}[G(\sm{X})P_m(\sm{X};\sm{\mathfrak{r}})].
\end{eqnarray*}
Again writing $\kappa_m$ for the coefficient of $\mathbb{E}[G(\sm{X})P_m(\sm{X};\sm{\mathfrak{r}})]$ in the above sum, we get that for any bounded measurable $G$,
$$\mathbb{E}[G(\sm{Y})] = \mathbb{E}\Big[G(\sm{X})\sum_m\kappa_mP_m(\sm{X};\sm{\mathfrak{r}})\Big].$$
Therefore $\sm{Y}\eqd\sm{X}$ if and only if $\sum_m\kappa_mP_m(\sm{X};\sm{\mathfrak{r}})=1$ or that for $m\geq1$,
$$\sum_{\sm{n}=m}H_n\frac1{\prod_{i=1}^dn_i!}\prod_{i=1}^d\mathfrak{h}_{n_i}(\mathfrak{r}_i)=0.$$

\end{enumerate}
\end{proof}

\begin{corollary}
If $d=2$ and $H_{n_1,n_2} = 0$ whenever $n_1\neq n_2$ so that expansion (\ref{expansion:0}) becomes
$$H(x_1,x_2) = \sum_{n\geq0}H_{n,n}P_n(x_1;\mathfrak{r}_1)P_n(x_2;\mathfrak{r}_2),$$
then $Y_1+Y_2$ has law $\mu(\mathfrak{r}_1)\ast\mu(\mathfrak{r}_2)$ if and only if $Y_1$ and $Y_2$ are independent.
\end{corollary}

\begin{corollary}
Suppose $Y_1$ and $Y_2$ have normal laws with parameters $\mathfrak{r}_1=(\mathfrak{m}_1,\mathfrak{s}_1^2)$ and $\mathfrak{r}_2=(\mathfrak{m}_2,\mathfrak{s}_2^2)$, respectively. Then $Y_1+Y_2$ is normal with parameter $\mathfrak{r}_1+\mathfrak{r}_2=(\mathfrak{m}_1+\mathfrak{m}_2,\mathfrak{s}_1^2+\mathfrak{s}_2^2)$ if and only if for every $n\geq1$,
$$\sum_{k=0}^nH_{k,n-k}\mathfrak{s}_1^{2k}\mathfrak{s}_2^{2(n-k)}=\sum_{k=1}^{n-1}H_{k,n-k}\mathfrak{s}_1^{2k}\mathfrak{s}_2^{2(n-k)}=0.$$
In particular, if $\mathfrak{s}_1=\mathfrak{s}_2=\mathfrak{s}$ then $Y_1+Y_2$ is normal with parameter $(\mathfrak{m}_1+\mathfrak{m}_2,2\mathfrak{s}^2)$ if and only if for every $n\geq1$,
$$\sum_{k=0}^nH_{k,n-k}=\sum_{k=1}^{n-1}H_{k,n-k}=0.$$
\end{corollary}
\begin{proof}
In this case $P_n(\cdot;\mathfrak{r})$ is the $n$th Hermite polynomial and $\mathfrak{h}_n(\mathfrak{r})=\mathfrak{s}^{2n}n!$. The result immediately follows by application of Theorem \ref{theorem:2}.
\end{proof}

For $x=(x_1,\ldots,x_d)$ and $A=\{i_1,\ldots,i_k\}\subset[d]$, we write $x_A$ for the reduced vector with elements indexed by $A$, $(x_{i_1},\ldots,x_{i_k})$. With this notation, if $Y=(Y_1,\ldots,Y_d)$ has law $Hd\mu_\otimes(\mathfrak{r})$ where $H$ is given by \eqref{expansion:0} and satisfies \eqref{expansion:0-assumption}, then $Y_A$ has law $H^{(A)}d\mu_\otimes(\mathfrak{r}_A)$ where
$H^{(A)}_{n_A} = H_{n_A^\circ}$ and $n_A^\circ$ has non-zero elements indexed by $A$ and zero elements indexed by $[d]\setminus A$ ($n_{\emptyset}^\circ=0$ and $n_{[d]}^\circ=n_{[d]}=n$).

\begin{corollary}
Let $A=\{i_1,\ldots,i_k\}\subset[d]$. Under the assumptions of Theorem \ref{theorem:2}, $\sm{Y_A}$ has law $\mu(\mathfrak{r}_{i_1})\ast\ldots\ast\mu(\mathfrak{r}_{i_k})$ if and only if for every $m\geq1$,
$$\sum_{\sm{n_A}=m}H_{n_A^\circ}\frac1{\prod_{i\in A}n_i!}\prod_{i\in A}\mathfrak{h}_{n_i}(\mathfrak{r}_i)=0.$$
\end{corollary}

\begin{remark}
Under the assumptions of Theorem \ref{theorem:2}, we know that $H_0=1$ and a necessary and sufficient condition for matching the marginal laws is that $H_n=0$ for any $n$ with a single non-zero element; that is if $n=n_A^\circ$ with $|A|=1$. Suppose further that $H_n=0$ whenever $n=n_A^\circ$ with $|A|>2$ (and $H\geq0$). Then $\sm{Y}\eqd\sm{X}$, where $X$ has law $\mu_\otimes(\mathfrak{r})$ and $Y$ has law $\mu_\bullet(\mathfrak{r})$, if and only if for every $m\geq1$,
$$\sum_{\substack{i,j\in[d]\\[1pt]i<j}}\sum_{n_i+n_j=m}H_{n_i,n_j}^{(\{i,j\})}\frac{\mathfrak{h}_{n_i}(\mathfrak{r}_i)\mathfrak{h}_{n_j}(\mathfrak{r}_j)}{n_i!n_j!}=0.$$
In particular, if $H_n=0$ whenever $n=n_A^\circ$ with $|A|>2$ (i.e. for any $n$ with more than 2 non-zero elements) and for every $m\geq1$ and every $i<j\leq d$,
$$\sum_{n_i+n_j=m}H_{n_i,n_j}^{(\{i,j\})}\frac{\mathfrak{h}_{n_i}(\mathfrak{r}_i)\mathfrak{h}_{n_j}(\mathfrak{r}_j)}{n_i!n_j!}=0$$
(and $H\geq0$), then for every $A\in[d]$, $\sm{Y_A}\eqd\sm{X_A}$, all sums of all subsets of $Y$ have laws that match those of $X$.

\end{remark}

\subsection{A construction using a finite expansion} Theorem \ref{theorem:2} gives necessary and sufficient conditions for $(Y_1,\ldots,Y_d)$ to match the marginal laws and the law of the sum of independent $(X_1,\ldots,X_d)$. This, however, does not lead to a construction of $Y=(Y_1,\ldots,Y_d)$ given $\mu_\otimes(\mathfrak{r})$, the law of $X=(X_1,\ldots,X_d)$. This is precisely what we attempt to do here. More specifically, given densities (within a given Meixner class) $\phi(\cdot;\mathfrak{r}_1),\ldots,\phi(\cdot;\mathfrak{r}_d)$, we look for a construction of a density of the form:
\begin{equation}\label{FinExpDensity}
\varphi(x) = \prod_{i=1}^d\phi(x_i;\mathfrak{r}_i) + \kappa\prod_{i=1}^d\overline{\phi}(x_i;\overline{\mathfrak{r}}_i)\sum_{1\leq\sm{n}\leq N}H_n\prod_{i=1}^d\overline{P}_{n_i}(x_i;\overline{\mathfrak{r}}_i),
\end{equation}
for quantities $N$, $(\kappa,H_n)$ (or simply $\kappa H_n$) and $(\overline{\phi},\overline{P}_n,\overline{\mathfrak{r}}_i)$ to be determined. Here $\overline{\phi}$ and $\overline{P}_n$ relate to a possibly different Meixner class $\mathcal{M}(\overline{\mathfrak{u}})$. The second term on the right is a perturbation of the density of $X$ and is chosen as to
\begin{enumerate}
\item integrate to 0 -- this is guaranteed by the orthogonality of the polynomials $P_{n_i}$;
\item be small enough for $\varphi$ to be nonnegative -- like in the Stoyanov example, $\kappa$ will play an important role here;
\item be such that $\sm{Y}\eqd\sm{X}$ -- condition \eqref{condition:0} will be crucial here.
\end{enumerate}

Next we state a variation of Theorem \ref{theorem:2}.
\begin{proposition}
Suppose that the function $\varphi$ defined in \eqref{FinExpDensity} is a density. Then the laws of the sums of $Y$ with density $\varphi(x)$ and $X$ with density $\prod_{i=1}^d\phi(x_i;\mathfrak{r}_i)$ are identical if and only if, for any $1\leq m\leq N$,
$$\sum_{\sm{n}=m}H_n\frac1{\prod_{i=1}^dn_i!}\prod_{i=1}^d\overline{\mathfrak{h}}_{n_i}(\overline{\mathfrak{r}}_i)=0.$$
\end{proposition}

\eqref{FinExpDensity} can be rewritten as
\begin{equation}\label{FinExpDensity2}
\varphi(x) = \prod_{i=1}^d\phi(x_i;\mathfrak{r}_i)\left(1 + \kappa\prod_{i=1}^d\frac{\overline{\phi}(x_i;\overline{\mathfrak{r}}_i)}{\phi(x_i;\mathfrak{r}_i)}\sum_{1\leq\sm{n}\leq N}H_n\prod_{i=1}^d\overline{P}_{n_i}(x_i;\overline{\mathfrak{r}}_i)\right)
\end{equation}
and $(\overline{\mathfrak{u}},\overline{\mathfrak{r}}_i)$ can be chosen (relative to $(\mathfrak{u},\mathfrak{r}_i)$) such that $|x_i|^N\overline{\phi}(x_i;\overline{\mathfrak{r}}_i)/\phi(x_i;\mathfrak{r}_i)\to 0$ as $|x_i|\to\infty$, then $\kappa$ can be chosen as to ensure that $\varphi$ is nonnegative. We discuss how this can be done for the first two Meixner types.

\begin{enumerate}

\item[(I)] \underline{$\mathfrak{a}=\mathfrak{b}=0$}. In this case $\phi$ is the density of a normal distribution and we choose $\overline{\mathfrak{u}}=\mathfrak{u}$ (no change in the Meixner class). Then, for $\mathfrak{r}=(\mathfrak{m},\mathfrak{s}^2)$, letting $\overline{\mathfrak{r}}=(\mathfrak{m},\overline{\mathfrak{s}}^2)$, where $\overline{\mathfrak{s}}^2<\mathfrak{s}^2$, leads to
$$|x|^N\frac{\overline{\phi}(x;\overline{\mathfrak{r}})}{\phi(x;\mathfrak{r})} = |x|^Ne^{-\frac12\left(1/\overline{\mathfrak{s}}^2-1/\mathfrak{s}^2\right)(x-\mathfrak{m})^2}\underset{|x|\to\infty}{\longrightarrow}0.$$

\item[(II)] \underline{$\mathfrak{a}=\mathfrak{b}<0$}. In this case $\phi$ is the density of a shifted gamma distribution
$$\phi(x;(\mathfrak{m},\mathfrak{s}^2)) = \frac{(-1/\mathfrak{a})^{\mathfrak{s}^2/\mathfrak{a}^2}}{\Gamma(\mathfrak{s}^2/\mathfrak{a}^2)}(x-\mathfrak{m}-\mathfrak{s}^2/\mathfrak{a})^{\mathfrak{s}^2/\mathfrak{a}^2-1}e^{x/\mathfrak{a}},\ x>\mathfrak{m}+\mathfrak{s}^2/\mathfrak{a};$$
$-\mathfrak{a}$ is the scale parameter, $\mathfrak{s}^2/\mathfrak{a}^2$ is the shape parameter and $\mathfrak{m}$ is a shift.
Here, we need to change Meixner classes at the same time as we change variances. We choose $\mathfrak{a}<\overline{\mathfrak{a}}<0$ (and $\overline{\mathfrak{b}}=\overline{\mathfrak{a}}$) and $\overline{\mathfrak{r}}$ such that $\overline{\mathfrak{s}}^2/\overline{\mathfrak{a}}^2=\mathfrak{s}^2/\mathfrak{a}^2$ with $\overline{\mathfrak{r}}^2<\mathfrak{s}^2$, and $\overline{\mathfrak{m}}+\overline{\mathfrak{s}}^2/\overline{\mathfrak{a}}=\mathfrak{m}+\mathfrak{s}^2/\mathfrak{a}$. We get
$$|x|^N\frac{\overline{\phi}(x;\overline{\mathfrak{r}})}{\phi(x;\mathfrak{r})} = |x|^N(\mathfrak{a}/\overline{\mathfrak{a}})^{\mathfrak{s}^2/\mathfrak{a}^2}e^{-x\left(1/\mathfrak{a}-1/\overline{\mathfrak{a}}\right)}
\underset{|x|\to\infty}{\longrightarrow}0.$$

\end{enumerate}

We now construct examples based on the normal distribution, Type (I), with zero mean (for simplicity). In this case Condition \eqref{condition:0} simplifies to
$$\sum_{\sm{n}=m}H_n\prod_{i=1}^d\mathfrak{s}_i^{2n_i}=0$$
and reduces further in the identically distributed case
$$\sum_{\sm{n}=m}H_n=0.$$

Analogous examples hold for generalized gamma distributions as well as other types within the Meixner family, such as the generalized binomial law ($\mathfrak{a}\neq\mathfrak{b}$, $\mathfrak{ab}<0$, $\mathfrak{a},\mathfrak{b}\in\mathbb{R}$, $\mathfrak{s}^2/(\mathfrak{a}\mathfrak{b})$ integer). In this case, the construction is even simpler because of the finiteness of the state space.

In the normal class the orthogonal polynomials $(P_n)_n$ are the Hermite polynomials, scaled to have unit leading coefficients. As the mean is fixed and equal to 0, we index these polynomials with the variance $\mathfrak{s}^2$ (instead of $\mathfrak{r}$).

\begin{example} \label{Example:2} Suppose $Y_1$ and $Y_2$ are marginally standard normal ($\mathfrak{r}=(0,1)$). The first three Hermite polynomials, scaled to have unit leading coefficients, are
$P_1(x;\overline{\mathfrak{s}}^2)= x,\> P_2(x;\overline{\mathfrak{s}}^2) =x^2-\overline{\mathfrak{s}}^2,\> P_3(x;\overline{\mathfrak{s}}^2) = x^3-3\overline{\mathfrak{s}}^2x$.
Choose $0<\overline{\mathfrak{s}}^2<1$ and, the non-zero terms to be $H_{1,3}=-1$ and $H_{3,1}=1$.
Then \eqref{FinExpDensity2} becomes
\begin{eqnarray}
\varphi(x_1,x_2) & = &\frac{e^{-\frac{1}{2}(x_1^2+x_2^2)}}{2\pi}\Big(1 + \kappa e^{-\frac12\left(1/\overline{\mathfrak{s}}^2-1\right)(x_1^2+x_2^2)}
\Big(-P_1(x_1;\overline{\mathfrak{s}}^2)P_3(x_2;\overline{\mathfrak{s}}^2)+P_3(x_1;\overline{\mathfrak{s}}^2)P_1(x_2;\overline{\mathfrak{s}}^2)\Big)\Big)
\nonumber \\
& = &\frac{e^{-\frac{1}{2}(x_1^2+x_2^2)}}{2\pi}\Big(1 + \kappa e^{-\frac12\left(1/\overline{\mathfrak{s}}^2-1\right)(x_1^2+x_2^2)}
x_1x_2\big(x_1^2-x_2^2\big)\Big),
\label{example:200}
\end{eqnarray}
which is equal to the density in the example by \citet{S2013} when $\overline{\mathfrak{s}}^2 = 1/2$.
\end{example}
\begin{example}\label{Example:3}
For an exchangeable solution consider non-zero terms $H_{1,3}=H_{3,1}=1$ and $H_{2,2}=-2$. Then
\begin{eqnarray*}
\sum_{\sm{n}=4}H_n\prod_{i=1}^2P_{n_i}(x_i;\overline{\mathfrak{s}}) & = & P_1(x_1;\overline{\mathfrak{s}}^2)P_3(x_2;\overline{\mathfrak{s}}^2)+P_3(x_1;\overline{\mathfrak{s}}^2)P_1(x_2;\overline{\mathfrak{s}}^2)-2P_2(x_1;\overline{\mathfrak{s}}^2)P_2(x_2;\overline{\mathfrak{s}}^2)\\
& = & x_1(x_2^3-3\overline{\mathfrak{s}}^2x_2)+x_2(x_1^3-3\overline{\mathfrak{s}}^2x_1)-2(x_1^2-\overline{\mathfrak{s}}^2)(x_2^2-\overline{\mathfrak{s}}^2)
\end{eqnarray*}
and
\begin{equation}
\varphi(x_1,x_2)=\frac{e^{-\frac{1}{2}(x_1^2+x_2^2)}}{2\pi}\Big(1 + \kappa e^{-\frac12\left(1/\overline{\mathfrak{s}}^2-1\right)(x_1^2+x_2^2)}
\big(x_1x_2(x_1^2+x_2^2-6\overline{\mathfrak{s}}^2) - 2(x_1^2-\overline{\mathfrak{s}}^2)(x_2^2-\overline{\mathfrak{s}}^2)\big)\Big).
\label{example:300}
\end{equation}	
Its moment generating function is,
\begin{equation}
M(t_1,t_2) = e^{(t_1^2+t_2^2)/2} + \kappa\overline{\mathfrak{s}}^8e^{\overline{\mathfrak{s}}^2(t_1^2+t_2^2)/2}t_1t_2(t_1-t_2)^2,
\label{trans:300}
\end{equation}
which follows from the transform of the $n$th Hermite polynomial
$$\mathbb{E}\big[e^{tX}P_n(X;\overline{\mathfrak{s}}^2)\big] = e^{\overline{\mathfrak{s}}^2t^2/2}\overline{\mathfrak{s}}^{2n}t^n.$$
As a check, the moment generating function of $Y_1+Y_2$ is correctly found to be $e^{-t^2}$ by setting $t_1=t_2=t$ in (\ref{trans:300}).
\end{example}
\begin{example}\label{Example:4}
Generalizing Example \ref{Example:3} under the assumption that $H_n=0$ for any $n$ with one or more than two non-zero elements or if $\sm{n}\neq4$, we let $Y$ be $d$-dimensional with density
\begin{eqnarray*}
\lefteqn{\varphi(x_1,\ldots,x_d)}\\
& = & \prod_{i=1}^d\phi(x_i;\mathfrak{s}_i^2) + \kappa\prod_{i=1}^d\phi(x_i;\overline{\mathfrak{s}}^2)\sum_{\substack{i,j\in [d]\\[1pt]i<j}}\sum_{n_i+n_j=4}H_{n_i,n_j}^{\{i,j\}}P_{n_i}(x_i;\overline{\mathfrak{s}}^2)P_{n_j}(x_j;\overline{\mathfrak{s}}^2)\\
& = & \prod_{i=1}^d\phi(x_i;\mathfrak{s}_i^2) + \kappa\prod_{i=1}^d\phi(x_i;\overline{\mathfrak{s}}^2)\sum_{\substack{i,j\in [d]\\[1pt]i<j}}\sum_{n_i+n_j=4}H_{n_i,n_j}^{\{i,j\}}\big(x_ix_j(x_i^2+x_j^2-6\overline{\mathfrak{s}}^2) - 2(x_i^2-\overline{\mathfrak{s}}^2)(x_j^2-\overline{\mathfrak{s}}^2)\big)
\end{eqnarray*}
where $0<\overline{\mathfrak{s}}^2<\min_{j\in [d]}\mathfrak{s}_j^2$ and $\kappa$ is sufficiently small. Suppose further that $H_{1,3}^{\{i,j\}}=H_{3,1}^{\{i,j\}}=1$ and $H_{2,2}^{\{i,j\}}=-2$ and
let $M_n(t;\mathfrak{s}^2) = \mathbb{E}\big[e^{tX}P_n(X;\mathfrak{s}^2)\big] = e^{\mathfrak{s}^2t^2/2}\mathfrak{s}^{2n}t^n$.
Then the moment generating function of $\varphi$ is
\begin{eqnarray*}
\lefteqn{M(t_1,\ldots,t_d)}\\
& = & \prod_{k=1}^de^{\mathfrak{s}_k^2t_k^2/2} + \kappa\sum_{\substack{i,j\in [d]\\[1pt]i<j}}\sum_{n_i+n_j=4}H_{n_i,n_j}^{\{i,j\}}M_{n_i}(t_i;\overline{\mathfrak{s}}^2)M_{n_j}(t_j;\overline{\mathfrak{s}}^2)\prod_{k\not\in\{i,j\}}e^{\overline{\mathfrak{s}}^2t_k^2/2}\\
& = & e^{(1/2)\sum_{k=1}^d\mathfrak{s}_k^2t_k^2}\\
& & +\ \kappa\sum_{\substack{i,j\in [d]\\[1pt]i<j}}\prod_{k\not\in\{i,j\}}e^{\overline{\mathfrak{s}}^2t_k^2/2}\Big(M_1(t_i;\overline{\mathfrak{s}}^2)M_3(t_j;\overline{\mathfrak{s}}^2)+M_3(t_i;\overline{\mathfrak{s}}^2)M_1(t_j;\overline{\mathfrak{s}}^2)-2M_2(t_i;\overline{\mathfrak{s}}^2)M_2(t_j;\overline{\mathfrak{s}}^2)\Big)\\
& = & e^{(1/2)\sum_{k=1}^d\mathfrak{s}_k^2t_k^2} + \kappa\overline{\mathfrak{s}}^8e^{(\overline{\mathfrak{s}}^2/2)\sum_{k=1}^dt_k^2}\sum_{\substack{i,j\in [d]\\[1pt]i<j}}t_it_j\big(t_i-t_j\big)^2.\\
\end{eqnarray*}
Then for any $A \subset[d]$, letting $t_i=t$ for $i\in A$ and, if $A\neq[d]$, $t_i=0$ for $i\notin A$ in the above, we see that $\sm{Y_A}\eqd\sm{Y_A}$. In other words all sub-sums of $Y$ have laws that match the same sub-sums for $X$. This is interesting as a comparison with a process with independent increments.	
\end{example}
\begin{remark}
Example \ref{Example:4} is suggestive of another approach to the construction of solutions, one based on moment generating functions of the form:
$$M(t_1,\ldots,t_d) = e^{(1/2)\sum_{k=1}^d\mathfrak{s}_k^2t_k^2} + \kappa e^{(\overline{\mathfrak{s}}^2/2)\sum_{k=1}^dt_k^2}\mathfrak{Q}(t_1,\ldots,t_d),$$
where $\mathfrak{Q}$ is a multivariate polynomial in $(t_1,\ldots,t_d)$. Using the identity, $\mathbb{E}\big[e^{tX}P_n(X;\mathfrak{s}^2)\big] = e^{\mathfrak{s}^2t^2/2}\mathfrak{s}^{2n}t^n$, it is then possible to express $\mathfrak{Q}(t_1,\ldots,t_d)$ as a transform of a multivariate polynomial $Q(x_1,\ldots,x_d)$, which in turn can be expanded in terms of Hermite polynomials. This formal inversion can be quite messy and is omitted here. Following the procedure adopted in Example \ref{Example:4}, one is then reduced to choosing $\kappa$ small enough to guarantee the positivity of $\varphi$, and to choosing $\mathfrak{Q}$ in such a way as to guarantee the matching of the marginal laws and the laws of the sums. In fact, this approach offers scope for considerably more general matchings. Indeed, consider the polynomial
$$\mathfrak{Q}(t_1,\ldots,t_d) = \sum_{k=1}^m\sum_{\substack{i,j\in A_k\\[1pt]i<j}}t^{(k)}_it^{(k)}_j\Big(t^{(k)}_i-t^{(k)}_j\Big)^2,$$
where $A_0,A_1,\ldots,A_m$ form a partition of $[d]$, each of $A_1,\ldots,A_m$ has at least two elements ($A_0$ may be empty), and the variables
$$t^{(0)}_1,\ldots,t^{(0)}_{|A_0|},t^{(1)}_1,\ldots,t^{(1)}_{|A_1|},t^{(2)}_1,\ldots,t^{(2)}_{|A_2|},\ldots,t^{(m)}_1,\ldots,t^{(m)}_{|A_m|}$$
form a rearrangement (permutation) of the variables $t_1,\ldots,t_d$.

Then the polynomial $\mathfrak{Q}$ vanishes if $t^{(k)}_i=0$ for every $k\in\{1,\ldots,m\}$ and $i\in\{1,\ldots,|A_k|\}$. It also vanishes if
for every $k\in\{1,\ldots,m\}$, $t^{(k)}_1=\ldots=t^{(k)}_{|A_k|}$. Letting $X^{(k)}_1,\ldots,X^{(k)}_{|A_k|}$ denote those random variables in $X$ that correspond to the variables $t^{(k)}_1,\ldots,t^{(k)}_{|A_k|}$, and similarly for $Y$, we see that, for $k\in\{1,\ldots,m\}$ and $i\in\{1,\ldots,|A_k|\}$,
$$(Y^{(0)}_1,\ldots,Y^{(0)}_{|A_0|},Y^{(k)}_i)\eqd(Y^{(0)}_1,\ldots,Y^{(0)}_{|A_0|},X^{(k)}_i)$$
and that
$$(Y^{(0)}_1,\ldots,Y^{(0)}_{|A_0|},\sm{Y_{A_1}},\ldots,\sm{Y_{A_m}})\eqd(X^{(0)}_1,\ldots,X^{(0)}_{|A_0|},\sm{X_{A_1}},\ldots,\sm{X_{A_m}}).$$
\end{remark}

\subsection{The (independent and) identically distributed case}

We end this section with a result that gives a necessary and sufficient condition for matching the laws of $g(Y)$ with those of $g(X)$, for any symmetric $g$, when $X_1,\ldots,X_d$ are independent as well as identically distributed random variables within one of the Meixner classes. The next section deals with the non-independent case by providing a generic construction.

We denote by $S_d$ the space of permutations on $[d]=\{1,\ldots,d\}$ and for $\beta\in S_d$, we write $\sigma_\beta$ for the function $\sigma_\beta(n)=(n_{\beta(1)},\ldots,n_{\beta(d)})$.

\begin{theorem}
Suppose $\mathfrak{r}_1=\ldots=\mathfrak{r}_d$ so that $\mathfrak{r}=(\mathfrak{r}_1,\ldots,\mathfrak{r}_1)$. Let $X$ have law $\mu_\otimes(\mathfrak{r})$. Then $Y$ with law $d\mu_\bullet(\mathfrak{r})=Hd\mu_\otimes(\mathfrak{r})$ satisfies $g(Y)\eqd g(X)$, for any $g$ symmetric, if and only if
\begin{equation}
\forall n\neq0,\ \sum_{\beta\in S_d}H_{\sigma_\beta(n)}=0.
\label{twocond:10}
\end{equation}
\end{theorem}
\begin{proof}
Suppose $g(Y)\eqd g(X)$, for any $g$ symmetric. Fix $n\neq0$ and let
$$g(x)=\sum_{\sigma\in S_d}\prod_{i=1}^dP_{n_{\sigma(i)}}(x_i;\mathfrak{r}_1).$$
Then $\mathbb{E}[g(X)]=0$ while
$$\mathbb{E}[g(Y)] = \mathbb{E}[g(X)H(X)] = \sum_{\beta\in S_d}\mathbb{E}\Big[H(X)\prod_{i=1}^dP_{n_{\beta(i)}}(X_i;\mathfrak{r}_1)\Big]
= \Big(\sum_{\beta\in S_d}H_{\sigma_\beta(n)}\Big)\prod_{i=1}^d\mathfrak{h}_{n_i}(\mathfrak{r}_1).$$
This trivially leads to \eqref{twocond:10}.
Conversely, if \eqref{twocond:10} holds, then for $G$ bounded, measurable and symmetric,
\begin{eqnarray*}
\lefteqn{\mathbb{E}[G(Y)]\ =\ \mathbb{E}[G(X)H(X)]}\\
& = & \frac1{d!}\sum_{\beta\in S_d}\mathbb{E}[G(X)H(X)]\ =\ \frac1{d!}\sum_{\beta\in S_d}\sum_nH_n\mathbb{E}\Big[G(X)\prod_{i=1}^dP_{n_i}(X_i;\mathfrak{r}_1)\Big]\\
& = & \frac1{d!}\sum_{\beta\in S_d}\sum_nH_{\sigma_\beta(n)}\mathbb{E}\Big[G(X)\prod_{i=1}^dP_{n_{\beta(i)}}(X_i;\mathfrak{r}_1)\Big]\\
& = & \frac1{d!}\sum_{\beta\in S_d}\sum_nH_{\sigma_\beta(n)}\mathbb{E}\Big[G(X)\prod_{i=1}^dP_{n_i}(X_{\beta^{-1}(i)};\mathfrak{r}_1)\Big]\\
& = & \frac1{d!}\sum_{\beta\in S_d}\sum_nH_{\sigma_\beta(n)}\mathbb{E}\Big[G(\sigma_{\beta^{-1}}(X))\prod_{i=1}^dP_{n_i}(X_{\beta^{-1}(i)};\mathfrak{r}_1)\Big]\\
& = & \frac1{d!}\sum_{\beta\in S_d}\sum_nH_{\sigma_\beta(n)}\mathbb{E}\Big[G(X)\prod_{i=1}^dP_{n_i}(X_i;\mathfrak{r}_1)\Big]\\
& = & \frac1{d!}\sum_n\Big(\sum_{\beta\in S_d}H_{\sigma_\beta(n)}\Big)\mathbb{E}\Big[G(X)\prod_{i=1}^dP_{n_i}(X_i;\mathfrak{r}_1)\Big].
\end{eqnarray*}
We end the proof by observing that $H_{\sigma_\beta(0)}=H_0=1$.
\end{proof}

\section{The identically distributed case -- A symmetry-balancing approach}\label{seccopula}

We shall assume throughout this section that all multivariate random variables admit joint densities, denoted by $f$.
We propose a generic construction under the assumption that the random variables $X_1,X_2,\ldots,X_d$ are identically distributed and admit a joint density but with no restriction on their dependence. An extension to non-indentically distributed random variables will be discussed at the end of the section (see Subsection \ref{nonidentical}). Although general, the proposed construction does not capture all possible solutions. A different approach is presented in Section \ref{secmeixner} where a full characterisation of the solution set is given.

\subsection{The two-dimensional case}\label{subsecdim2}

Let $X_1$ and $X_2$ be two identically distributed random variables with joint density $f$, marginal distribution function $\Phi$, assumed to be continuous and strictly increasing, and marginal density $\phi$. Let $\mathcal{G}_2$ be the set of symmetric functions on $\mathbb{R}^2$, that is functions $g$ that satisfy the property that
$$\forall(x_1,x_2)\in\mathbb{R}^2,\ g(x_2,x_1)=g(x_1,x_2),$$
and let
$$D_\pm = \{(x_1,x_2):\pm x_2<\pm x_1\}.$$
Then, for any (measurable) $\ell$ such that $0\leq\ell\leq f$,
$$\varphi(x_1,x_2) =
\begin{cases}
f(x_1,x_2)-\ell(x_1,x_2) & (x_1,x_2)\in D_+\\
f(x_1,x_2)+\ell(x_2,x_1) & (x_1,x_2)\in D_-\\
\end{cases}$$
is the density of a pair $(Y_1,Y_2)$ such that for any $g\in\mathcal{G}_2$, $g(Y_1,Y_2)\eqd g(X_1,X_2)$.

While $\varphi$ defines a pair that that matches the law of $g(X_1,X_2)$, for all symmetric functions $g$, it does not necessarily preserve the marginal laws of $f$. To do so we need to also ``compensate'' in the $x_1$ and $x_2$ directions. We illustrate this point by thinking of the density of the uniform distribution on the unit square perturbed with a hump (in the $x_1<x_2$ region) and a corresponding trough (in the $x_1>x_2$ region). While such a perturbation maintains say the law of the sum, it alters marginal laws. To maintain those, we must compensate the troughs and humps in both the $x_1$ and $x_2$ directions -- See Figure \ref{figNoMarg} below. We call this a symmetry-balancing approach.

\begin{figure}[!htb]
\begin{center}
\includegraphics[height=6cm]{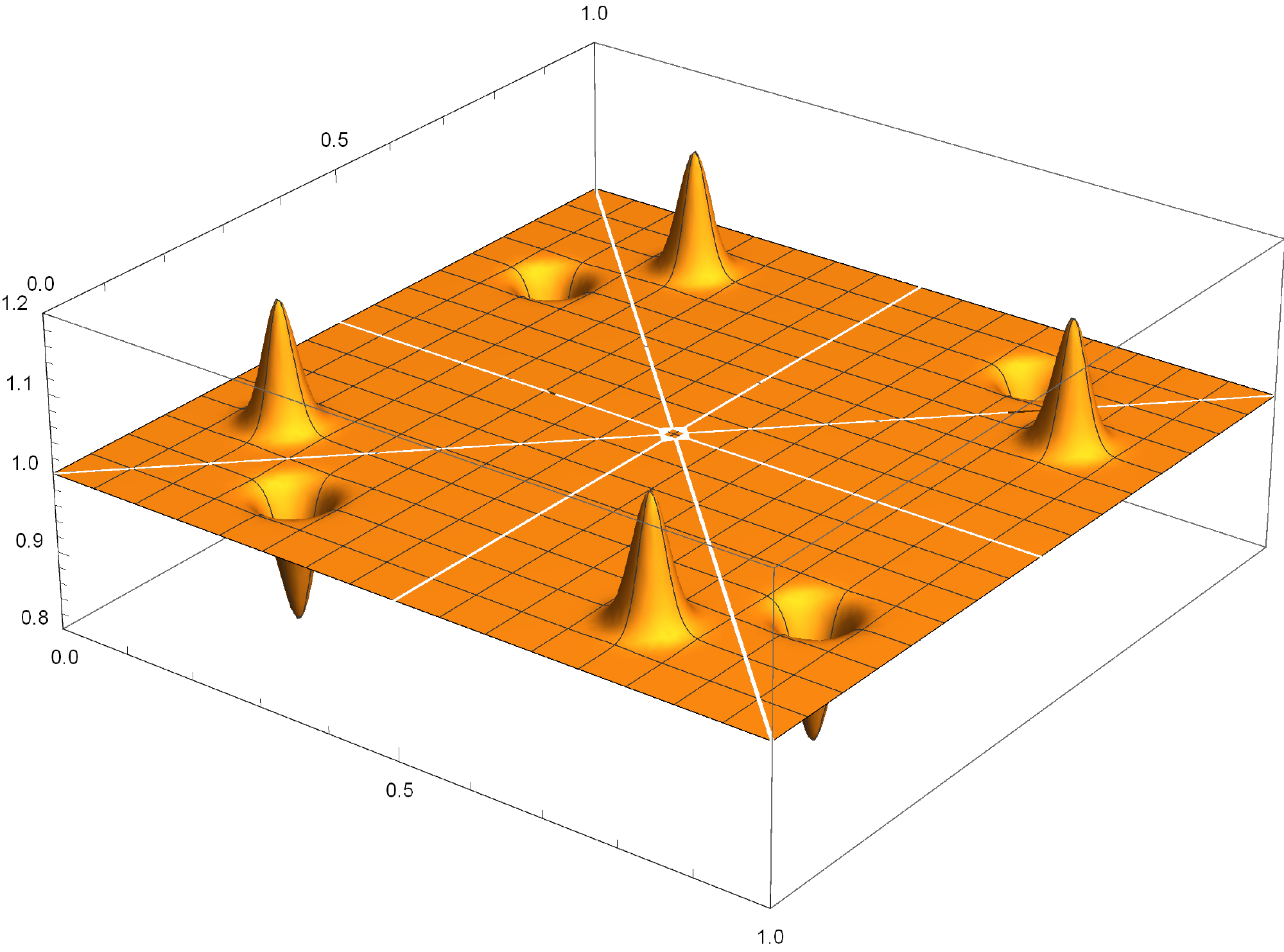}
\end{center}
\caption{The symmetry-balancing approach in the case of the uniform distribution on the unit square}
\label{figNoMarg}
\end{figure}

The proposed construction is universal in that it is described through the use of a ``copula perturbation'' that can then be applied to any distribution. Let $c$ be the copula density of $f$.

\begin{proposition}\label{propOctalCopula}
For any (measurable) $\gamma:[0,1]\times[0,1]\to[0,+\infty)$ with
\begin{equation}\label{gammacond0}
\gamma(u_1,u_2)\leq \min(c(u_1,u_2),c(u_2,1-u_1),c(1-u_2,u_1),c(1-u_1,1-u_2))
\end{equation}
the function
\begin{equation}
\theta_\gamma(u_1,u_2) = \begin{cases}
c(u_1,u_2)-\gamma(u_1,u_2) & u_1\leq u_2\leq \frac{1}{2} \\
c(u_1,u_2)+\gamma(u_1,1-u_2) & \frac{1}{2}\leq u_2\leq 1-u_1 \\
c(u_1,u_2)-\gamma(1-u_2,u_1) & 1-u_2\leq u_1\leq \frac{1}{2} \\
c(u_1,u_2)+\gamma(1-u_2,1-u_1) & \frac{1}{2}\leq u_1\leq u_2 \\
c(u_1,u_2)-\gamma(1-u_1,1-u_2) & \frac{1}{2}\leq u_2\leq u_1 \\
c(u_1,u_2)+\gamma(1-u_1,u_2) & 1-u_1\leq u_2\leq \frac{1}{2} \\
c(u_1,u_2)-\gamma(u_2,1-u_1) & \frac{1}{2}\leq u_1\leq 1-u_2 \\
c(u_1,u_2)+\gamma(u_2,u_1) & u_2\leq u_1\leq \frac{1}{2}
\end{cases}
\label{eqOctalCopDef}
\end{equation}
is a copula density.
Furthermore, for any bounded and symmetric $g$,
$$\int g(u_1,u_2)\theta_\gamma(u_1,u_2)du_1dv = \int g(u_1,u_2)c(u_1,u_2)du_1dv.$$
\end{proposition}
The proof is given in Theorem \ref{thmGenCopula} in the more general setting of the $d$-dimensional case, $d\geq2$.

When $c$ is the density of the independence copula ($c=1$), we call the copula associated with $\theta_\gamma$, the ``Octal Copula with generator $\gamma$''.

Finally, we note that the Stoyanov example has a copula that is of the octal form.

\begin{remark}
Let $Y_1$ and $Y_2$ be identically distributed and belong to one of the Meixner classes. Suppose that they are jointly absolutely continuous and admit an octal copula density constructed from \eqref{eqOctalCopDef} by taking $c(u,v)=1$.
By virtue of the properties of the octal copula, for $n_1+n_2\geq1$,
$$(H_{n_1,n_2}+H_{n_2,n_1})\mathfrak{h}_{n_1}(\mathfrak{r})\mathfrak{h}_{n_2}(\mathfrak{r})=\mathbb{E}\big [P_{n_1}(Y_1;\mathfrak{r})P_{n_2}(Y_2;\mathfrak{r}) +P_{n_1}(Y_2;\mathfrak{r})P_{n_2}(Y_1;\mathfrak{r})\big ]=0,$$
as the expectation of a symmetric function of $(Y_1,Y_2)$ which has the same expected value 0 as if $Y_1$ and $Y_2$ were independent. In other words (\ref{twocond:10}) holds true.

Next we assume that $Y_1$ and $Y_2$ are marginally standard normal and proceed to identify the coefficients $H_{n_1,n_2}$ in this case. Here $\mathfrak{r}=(0,1)$ and we shall simply write $P_n$ for $P_n(\cdot;\mathfrak{r})$ and $\mathfrak{h}_n$ for $\mathfrak{h}_n(\mathfrak{r})$. We know that $\Phi(-x) = 1 - \Phi(x)$, $P_n(-x) = (-1)^nP_n(x)$, a property of the Hermite polynomials that can easily be checked. Recall that
\begin{eqnarray*}
\theta_\gamma(u_1,u_2)+\theta_\gamma(u_2,u_1) & = & \theta_\gamma(u_1,u_2)+\theta_\gamma(1-u_1,u_2)\ =\ \theta_\gamma(u_1,u_2)+\theta_\gamma(u_1,1-u_2)\\
& = & \theta_\gamma(u_1,u_2)+\theta_\gamma(1-u_2,1-u_1)\ =\ 2
\end{eqnarray*}
and that
$$\theta_\gamma(u_1,u_2)-\theta_\gamma(1-u_2,u_1)=\theta_\gamma(u_1,u_2)-\theta_\gamma(u_2,1-u_1)=\theta_\gamma(u_1,u_2)-\theta_\gamma(1-u_1,1-u_2)=2.$$
These identities translates in terms of $H_{n_1,n_2}$, for $n_1+n_2\geq1$, to
$$H_{n_1,n_2}+H_{n_2,n_1}=\big(1+(-1)^{n_1}\big)H_{n_1,n_2}=\big(1+(-1)^{n_2}\big)H_{n_1,n_2}=H_{n_1,n_2}+(-1)^{n_1+n_2}H_{n_2,n_1}=0$$
and that
$$H_{n_1,n_2}-(-1)^{n_1}H_{n_2,n_1}=H_{n_1,n_2}-(-1)^{n_2}H_{n_2,n_1}=\big(1-(-1)^{n_1+n_2}\big)H_{n_1,n_2}=0.$$
We immediately deduce that $H_{n_1,n_2}=0$ if either $n_1$ is even or $n_2$ is even, and $H_{n_1,n_2}+H_{n_2,n_1}=0$ when both are odd.
\end{remark}

\subsection{The multidimensional ``matching'' copula}

We are now ready to extend the construction of $\theta_\gamma$ to the $d$-dimensional hypercube $[0,1]^d$ and consequently $\varphi$ to the $d$-dimensional Euclidean space $\mathbb{R}^d$. We shall retain from the two-dimensional case the idea that regions ($\Delta_2,\ldots,\Delta_8$) are mapped onto a reference region ($\Delta_1$). Core to these mappings are the reflections $(u_1,u_2)\hookrightarrow(1-u_1,u_2)$, $(u_1,u_2)\hookrightarrow(u_1,1-u_2)$ and $(u_1,u_2)\hookrightarrow(1-u_1,1-u_2)$ as well as $(u_1,u_2)\hookrightarrow(u_2,u_1)$. Generalising these to the hypercube lead to the maps $\tau_\alpha$, for the first three, and the maps $\sigma_\beta$ for the last one. These are introduced next.

\begin{itemize}
\item Recall that $S_d$ is the space of permutations on $[d]=\{1,\ldots,d\}$, and that for $\beta\in S_d$, $\sigma_\beta$ denotes the function
$$\sigma_\beta(u)=(u_{\beta(1)},\ldots,u_{\beta(d)}).$$
\item $id$ denotes the identity function.
\item $\mathcal{G}(\mathbb{R}^d)$ and $\mathcal{G}([0,1]^d)$, or simply $\mathcal{G}_d$, denote the sets of symmetric functions $g$ on $\mathbb{R}^d$ and $[0,1]^d$, respectively; that is for any $\beta\in S_d$, $g\circ\sigma_\beta=g$.
\item For $\alpha\in\{0,1\}^d$ and $u=(u_1,\ldots,u_d)\in[0,1]^d$,
$$\tau_\alpha(u) = \big(\alpha_i(1-u_i)+(1-\alpha_i)u_i\big)_{i=1}^d.$$
\item $\Delta(0) = (0,1/2)^d$ and for any other $\alpha\in\{0,1\}^d$,
$$\Delta(\alpha) = \{u\in[0,1]^d:\ \tau_\alpha(u)\in\Delta(0)\}.$$
\item $\Delta(0,id) = \{u\in[0,1]^d:\ u_1<u_2<\ldots<u_d<1/2\}$ and, for any other pair $(\alpha,\beta)\in\{0,1\}^d\times S_d$,
$$\Delta(\alpha,\beta) = \{u\in[0,1]^d:\ \sigma_\beta(\tau_\alpha(u))\in\Delta(0,id)\}.$$
In the two-dimensional case of Subsection \ref{subsecdim2}, $\Delta(0,id)$ was referred to as $\Delta_1$, $\Delta(0,\perm{21})$ was referred to as $\Delta_8$, $\Delta((0,1),id)$ was referred to as $\Delta_2$ etc.
\item $\Xi_d=\{u\in[0,1]^d:\exists k,k'\in[d]\mbox{ such that }(u_k-1/2)(u_k-u_{k'})=0\}$.
\end{itemize}

The next lemma shows that we can essentially partition the hypercube into $2^dd!$ regions that all map onto the reference region $\Delta(0,id)$.
\begin{lemma}\label{lemmain}
\begin{enumerate}
\item $\bigcup_{(\alpha,\beta)}\Delta(\alpha,\beta) = [0,1]^d\setminus\Xi_d$; i.e. $\forall u\in[0,1]^d\setminus\Xi_d$, $\exists!(\alpha,\beta)\in\{0,1\}^d\times S_d$ such that $q=\sigma_\beta(\tau_\alpha(u))$ satisfies the condition
$0<q_1<\ldots<q_d<1/2$.
\item For $(\alpha,\beta)\neq(\alpha',\beta')$, $\Delta(\alpha,\beta)\cap\Delta(\alpha',\beta')=\emptyset$.
\item For any $u\in[0,1]^d$, $\displaystyle\sigma_\beta(\tau_\alpha(u)) = \tau_{\sigma_\beta(\alpha)}(\sigma_\beta(u))$.
\end{enumerate}
\end{lemma}
\begin{proof}
(1) and (2) are immediate.

For any $u\in[0,1]^d$,
\begin{eqnarray*}
\sigma_\beta(\tau_\alpha(u)) & = & \sigma_\beta\Big(\big(\alpha_i(1-u_i)+(1-\alpha_i)u_i\big)_{i=1}^d\Big)\\
& = & \Big(\alpha_{\beta(i)}(1-u_{\beta(i)})+(1-\alpha_{\beta(i)})u_{\beta(i)}\Big)_{i=1}^d\ =\ \tau_{\sigma_\beta(\alpha)}(\sigma_\beta(u)).
\end{eqnarray*}
This proves (3).
\end{proof}

A necessary step in the construction is the embedding of the $(d-1)$-dimensional hypercube as a hyperplane in the $d$-dimensional hypercube and how the corresponding partitions carry across.
To that end, we need the following notations and results.
\begin{itemize}
\item For $v\in[0,1]^{d-1}$, $r\in[0,1]$ and $k\in[d]$,
$$\omega_k(v,r)=(v_1,\ldots,v_{k-1},r,v_k,\ldots,v_{d-1}),$$
with the obvious adjustments for the cases $k=1$ and $k=d$.
\item For $b\in S_{d-1}$ and $j,k\in[d]$, $\beta=\chi_j(b,k)$ is the permutation in $S_d$
$$\beta(i)=\left\{
\begin{array}{ll}
b(i) & \text{if } i\leq j-1\mbox{ and }b(i)\leq k-1\\
b(i)+1 & \text{if } i\leq j-1\mbox{ and }b(i)\geq k\\
k & \text{if } i=j\\
b(i-1) & \text{if } i\geq j+1\mbox{ and }b(i-1)\leq k-1\\
b(i-1)+1 & \text{if } i\geq j+1\mbox{ and }b(i-1)\geq k
\end{array}
\right.$$
\item For $a\in\{0,1\}^{d-1}$, $b\in S_{d-1}$, $j\in[d]$ and $r\in[0,1]$, with a slight abuse of notation in the use of $\tau$ and $\sigma$,
$$\Delta_j(a,b,r) = \{v\in[0,1]^{d-1}:\ q_0<\ldots<q_{j-1}<\bar r<q_j<\ldots<q_d\},$$
where $\bar r=\min(r,1-r)$ and, $q$ stands for $\sigma_b(\tau_a(v))$ and has been augmented with the bounds $q_0=0$ and $q_d=1/2$.
\end{itemize}

\begin{lemma}\label{lemmarginal}
Let $a\in\{0,1\}^{d-1}$, $b\in S_{d-1}$ and $j,k\in[d]$.
\begin{enumerate}
\item If $\beta=\chi_j(b,k)$ and $r\in[0,1]$, then
$$\tau_\alpha(\omega_k(v,r)) = \omega_k(\tau_a(v),\bar r)\mbox{ and }\sigma_\beta(\tau_\alpha(\omega_k(v,r))) = \omega_j(\sigma_b(\tau_a(v)),\bar r),$$
where
$$\alpha=\left\{
\begin{array}{ll}
\omega_k(a,0) & r\in[0,1/2]\\
\omega_k(a,1) & r\in(1/2,1]\\
\end{array}\right.$$
\item With $\alpha$ and $\beta$ as above,
$$\omega_k(v,r)\in\Delta(\alpha,\beta)\quad\Longleftrightarrow\quad v\in\Delta_j(a,b,r)\quad\Longleftrightarrow\quad\sigma_b(\tau_a(v))\in\Delta_j(0,id,r).$$
\end{enumerate}
\end{lemma}

Since $\Xi_d$ has Lebesgue measure zero, any integral over $[0,1]^d$ can be taken to mean an integral on $[0,1]^d\setminus\Xi_d$.

\begin{theorem}\label{thmGenCopula}
Let $c$ be a density on $[0,1]^d$, $\gamma$ be an integrable non-negative function on $\Delta(0,id)$ and $\varepsilon:\{0,1\}^d\times S_d\to[-1,1]$. We assume that
$$\gamma(u)\leq \min_{\varepsilon(\alpha,\beta)>0}\frac{c(\tau_\alpha(\sigma_{\beta^{-1}}(u)))}{\varepsilon(\alpha,\beta)}.$$
Define on $[0,1]^d$ the function
\begin{equation}
\theta_{\varepsilon,\gamma}(u) = c(u)-\sum_{\alpha,\beta}\varepsilon(\alpha,\beta)\gamma(\sigma_\beta(\tau_\alpha(u)))1_{\Delta(\alpha,\beta)}(u).
\label{eqGenCopDef}
\end{equation}
\begin{enumerate}
\item[(D)] $\theta_{\varepsilon,\gamma}$ is a density if and only if
\begin{equation}\label{epscondD}
\sum_{\alpha,\beta}\varepsilon(\alpha,\beta)=0.
\end{equation}
\end{enumerate}
Assume \eqref{epscondD} and let $U$ and $V$ be multivariate random variables with densities $c$ and $\theta_{\varepsilon,\gamma}$, respectively.
\begin{enumerate}
\item[(G)] $g(V)\eqd g(U)$, for any $g$ symmetric, if and only if
\begin{equation}\label{epscondG}
\forall\alpha\in\{0,1\}^d,\ \sum_{\beta}\varepsilon(\sigma_{\beta^{-1}}(\alpha),\beta)=0.
\end{equation}
In particular, in this case, $V_1+\ldots+V_d\eqd U_1+\ldots+U_d$.

\item[($\mathfrak{M}$)] Fix $k\in[d]$. If
\begin{equation}\label{epscondM}
\forall j\in[d],\forall a\in\{0,1\}^{d-1},\forall b\in S_{d-1},\ \sum_{\aleph=0,1}\varepsilon(\omega_k(a,\aleph),\chi_j(b,k))=0,
\end{equation}
then $(V_1,\ldots,V_{k-1},V_{k+1},\ldots,V_d) \eqd (U_1,\ldots,U_{k-1},U_{k+1},\ldots,U_d)$.
\item[($\mathfrak{C}$)] Suppose $c$ is a copula density. If
\begin{equation}\label{epscondC}
\forall j,k\in[d],\ \sum_{a,b}\varepsilon(\omega_k(a,0),\chi_j(b,k))=\sum_{a,b}\varepsilon(\omega_k(a,1),\chi_j(b,k))=0,
\end{equation}
then $\theta_{\varepsilon,\gamma}$ is a copula density.
\end{enumerate}
\end{theorem}

\begin{proof}
\begin{enumerate}
\item[(D)] Clearly, $\theta_{\varepsilon,\gamma}(u)\geq0$. Let us show that it integrates to 1 or equivalently that \linebreak $\sum_{\alpha,\beta}\varepsilon(\alpha,\beta)\gamma(\sigma_\beta(\tau_\alpha(u)))1_{\Delta(\alpha,\beta)}(u)$ integrates to 0.

We show in (G) below that for any symmetric $g$,
\begin{eqnarray*}
\lefteqn{\int_{[0,1]^d}g(u)\theta_{\varepsilon,\gamma}(u)du}\\
& = & \int_{[0,1]^d}g(u)c(u)du-\sum_\alpha\int_{\Delta(0,id)}g(\tau_{\alpha}(u))\gamma(u)du\sum_\beta\varepsilon(\sigma_{\beta^{-1}}(\alpha),\beta).
\end{eqnarray*}
Applied to $g=1$, we get that $\theta_{\varepsilon,\gamma}$ is a density if and only if
$$0=\sum_\alpha\sum_\beta\varepsilon(\sigma_{\beta^{-1}}(\alpha),\beta)=\sum_\beta\sum_\alpha\varepsilon(\sigma_{\beta^{-1}}(\alpha),\beta)=\sum_{\alpha,\beta}\varepsilon(\alpha,\beta).$$

\item[(G)] For any symmetric function $g$,
\begin{eqnarray*}
\lefteqn{\int_{[0,1]^d}g(u)\theta_{\varepsilon,\gamma}(u)du}\\
& = & \int_{[0,1]^d}g(u)c(u)du-\sum_{\alpha,\beta}\varepsilon(\alpha,\beta)\int_{\Delta(\alpha,\beta)}g(u)\gamma(\sigma_\beta(\tau_\alpha(u)))du\\
& = & \int_{[0,1]^d}g(u)c(u)du-\sum_{\alpha,\beta}\varepsilon(\alpha,\beta)\int_{\Delta(0,id)}g(\tau_\alpha(\sigma_{\beta^{-1}}(u)))\gamma(u)du\\
& = & \int_{[0,1]^d}g(u)c(u)du-\sum_{\alpha,\beta}\varepsilon(\alpha,\beta)\int_{\Delta(0,id)}g(\sigma_{\beta^{-1}}(\tau_{\sigma_\beta(\alpha)}(u)))\gamma(u)du\\
& = & \int_{[0,1]^d}g(u)c(u)du-\sum_{\alpha,\beta}\varepsilon(\alpha,\beta)\int_{\Delta(0,id)}g(\tau_{\sigma_\beta(\alpha)}(u))\gamma(u)du\\
& = & \int_{[0,1]^d}g(u)c(u)du-\sum_{\alpha'}\int_{\Delta(0,id)}g(\tau_{\alpha'}(u))\gamma(u)du\sum_{\alpha,\beta:\sigma_\beta(\alpha)=\alpha'}\varepsilon(\alpha,\beta)\\
& = & \int_{[0,1]^d}g(u)c(u)du-\sum_\alpha\int_{\Delta(0,id)}g(\tau_{\alpha}(u))\gamma(u)du\sum_\beta\varepsilon(\sigma_{\beta^{-1}}(\alpha),\beta).
\end{eqnarray*}
Clearly if \eqref{epscondG} holds, then $\displaystyle\int_{[0,1]^d}g(u)\theta_{\varepsilon,\gamma}(u)du=\int_{[0,1]^d}g(u)c(u)du$ and $g(V)\eqd g(U)$. As already observed, the latter point follows from the fact that if $G$ is bounded and $g$ is symmetric, then $G\circ g$ is bounded and symmetric.

Conversely, suppose that $g(V)\eqd g(U)$, for any $g$ symmetric; i.e. for any $g$ symmetric,
$$\sum_\alpha\int_{\Delta(\alpha,id)}\Gamma(\alpha)g(u)\gamma(\tau_\alpha(u))du=0,$$
where $\Gamma(\alpha)=\sum_\beta\varepsilon(\sigma_{\beta^{-1}}(\alpha),\beta)$. Fix $\alpha_0\in\{0,1\}^d$ and let $e_0$ be any bounded measurable function. Define $e$ as
$$e(u) = \left\{
\begin{array}{ll}
e_0(u) & u\in\Delta(\alpha_0,id)\\
0 & u\not\in\Delta(\alpha_0,id)
\end{array}
\right.$$
and $g$ by symmetrisation of $e$
$$g(u)=\sum_\beta e(\sigma_\beta(u)).$$
Then ($g$ is symmetric and)
\begin{eqnarray*}
0 & = & \sum_\alpha\int_{\Delta(\alpha,id)}\Gamma(\alpha)g(u)\gamma(\tau_\alpha(u))du\\
& = & \sum_\alpha\int_{\Delta(\alpha,id)}\Gamma(\alpha)\Big(\sum_\beta e(\sigma_\beta(u))\Big)\gamma(\tau_\alpha(u))du\\
& = & \sum_\alpha\sum_\beta\int_{\Delta(\alpha,id)}\Gamma(\alpha)e(\sigma_\beta(u))\gamma(\tau_\alpha(u))du\\
& = & \sum_\alpha\sum_\beta\int_{\Delta(\alpha,id)\cap\Delta(\sigma_{\beta^{-1}}(\alpha_0),\beta)}\Gamma(\alpha_0)e(\sigma_\beta(u))\gamma(\tau_{\alpha_0}(u))du\\
& = & \Gamma(\alpha_0)\sum_\beta\int_{\Delta(\alpha_0,id)}e(\sigma_\beta(u))\gamma(\tau_{\alpha_0}(u))du,
\end{eqnarray*}
where the second last identity follows from the fact that $\sigma_\beta(u)\in\Delta(\alpha_0,id)$ if and only if $u\in\Delta(\sigma_{\beta^{-1}}(\alpha_0),\beta)$.
We immediately conclude that $\Gamma(\alpha_0)=0$. Repeating for all other $\alpha_0$ in $\{0,1\}^d$, we prove the result.

\item[($\mathfrak{M}$)] Fix $k\in[d]$, $v\in[0,1]^{d-1}\setminus\Xi_{d-1}$ and let $(a,b)\in\{0,1\}^{d-1}\times S_{d-1}$ such that $q=\sigma_b(\tau_a(v))$ satisfies $0<q_1<\ldots<q_{d-1}<1/2$. Let $q_0=0$, $q_d=1/2$.

Using Lemma \ref{lemmarginal}, we get
\begin{eqnarray*}
\lefteqn{\int_{[0,1]}\theta_{\varepsilon,\gamma}(\omega_k(v,r))dr\ =\ \int_{[0,1]}c(\omega_k(v,r))dr}\\
& & -\ \sum_{j=0}^{d-1}\int_{[q_j,q_{j+1}]}\sum_{\alpha,\beta}\varepsilon(\alpha,\beta)\gamma(\sigma_\beta(\tau_\alpha(\omega_k(v,r))))1_{\Delta(\alpha,\beta)}(\omega_k(v,r))dr\\
& & -\ \sum_{j=0}^{d-1}\int_{[1-q_{j+1},1-q_j]}\sum_{\alpha,\beta}\varepsilon(\alpha,\beta)\gamma(\sigma_\beta(\tau_\alpha(\omega_k(v,r))))1_{\Delta(\alpha,\beta)}(\omega_k(v,r))dr\\
& = & \int_{[0,1]}c(\omega_k(v,r))dr\\
& & -\ \sum_{j=0}^{d-1}\int_{[q_j,q_{j+1}]}\varepsilon(\omega_k(a,0),\chi_{j+1}(b,k))\gamma(\omega_{j+1}(\sigma_b(\tau_a(v)),r))dr\\
& & -\ \sum_{j=0}^{d-1}\int_{[1-q_{j+1},1-q_j]}\varepsilon(\omega_k(a,1),\chi_{j+1}(b,k))\gamma(\omega_{j+1}(\sigma_b(\tau_a(v)),1-r))dr\\
& = & \int_{[0,1]}c(\omega_k(v,r))dr\\
& & -\ \sum_{j=0}^{d-1}\int_{[q_j,q_{j+1}]}\gamma(\omega_{j+1}(\sigma_b(\tau_a(v)),r))dr\sum_{\aleph=0,1}\varepsilon(\omega_k(a,\aleph),\chi_{j+1}(b,k))\\
& & \qquad
\end{eqnarray*}
If \eqref{epscondM} holds, then for any $k\in[d]$ and $v\in[0,1]^{d-1}\setminus\Xi_{d-1}$,
$$\int_{[0,1]}\theta_{\varepsilon,\gamma}(\omega_k(v,r))dr=\int_{[0,1]}c(\omega_k(v,r))dr;$$
that is $(V_1,\ldots,V_{k-1},V_{k+1},\ldots,V_d)\eqd(U_1,\ldots,U_{k-1},U_{k+1},\ldots,U_d)$.

%Conversely, suppose that for any $v\in[0,1]^{d-1}$,
%$$\sum_{j=0}^{d-1}\int_{[q_j,q_{j+1}]}\gamma(\omega_{j+1}(\sigma_b(\tau_a(v)),r))dr\sum_{\eta=0,1}\varepsilon(\omega_k(a,\eta),\chi_{j+1}(b,k))=0.$$
%Let $q\to\vartheta_j=v_j=\mbox{v}_j=\mathbf{v}_j$, the $j$th vertex ($\vartheta_1=\ldots=\vartheta_j=0$ and $\vartheta_{j+1}=\ldots=\vartheta_{d-1}=1/2$).

\item[($\mathfrak{C}$)] We fix $k\in[d]$ and proceed to show that under \eqref{epscondC}, $V_k\eqd U_k$. Fix $r\in(0,1)$.
\begin{eqnarray*}
\lefteqn{\int_{[0,1]^{d-1}}\theta_{\varepsilon,\gamma}(\omega_k(v,r))dv}\\
& = & \int_{[0,1]^{d-1}}c(\omega_k(v,r))dv\\
& & -\ \int_{[0,1]^{d-1}}\sum_{\alpha,\beta}\varepsilon(\alpha,\beta)\gamma(\sigma_\beta(\tau_\alpha(\omega_k(v,r))))1_{\Delta(\alpha,\beta)}(\omega_k(v,r))dv\\
& = & 1-\int_{[0,1]^{d-1}}\sum_{j,a,b}\varepsilon(\alpha,\beta)\gamma(\sigma_\beta(\tau_\alpha(\omega_k(v,r))))1_{\Delta_j(a,b,r)}(v)dv,
\end{eqnarray*}
where as in Lemma \ref{lemmarginal},
$$\alpha=\left\{
\begin{array}{ll}
\omega_k(a,0) & r\in[0,1/2]\\
\omega_k(a,1) & r\in(1/2,1]\\
\end{array}\right.
\mbox{ and }\beta=\chi_j(b,k).$$
It follows that
\begin{eqnarray*}
\int_{[0,1]^{d-1}}\theta_{\varepsilon,\gamma}(\omega_k(v,r))dv & = & 1-\sum_{j,a,b}\varepsilon(\alpha,\beta)\int_{\Delta_j(a,b,r)}\gamma(\sigma_\beta(\tau_\alpha(\omega_k(v,r))))dv\\
& = & 1-\sum_{j,a,b}\varepsilon(\alpha,\beta)\int_{\Delta_j(a,b,r)}\gamma(\omega_j(\sigma_b(\tau_a(v)),\bar r))dv\\
& = & 1-\sum_{j,a,b}\varepsilon(\alpha,\beta)\int_{\Delta_j(0,id,r)}\gamma(\omega_j(v,\bar r))dv\\
& = & 1-\sum_{j=1}^d\int_{\Delta_j(0,id,r)}\gamma(\omega_j(v,\bar r))dv\sum_{a,b}\varepsilon(\alpha,\beta).
\end{eqnarray*}
which concludes the proof.

\end{enumerate}
\end{proof}

\begin{example}[The bivariate case]
When $d=2$, $S_2$ consists of the identity permutation and the transposition $\perm{21}$. In this case, \eqref{epscondC} produces 2 sets of 4 equations:
$$\mbox{for }k=1\left\{
\begin{array}{l}
\varepsilon((0,0),\perm{21})+\varepsilon((0,1),\perm{21}) = 0\\
\varepsilon((0,0),id)+\varepsilon((0,1),id) = 0\\
\varepsilon((1,0),\perm{21})+\varepsilon((1,1),\perm{21}) = 0\\
\varepsilon((1,0),id)+\varepsilon((1,1),id) = 0
\end{array}
\right.\mbox{ and for }k=2\left\{
\begin{array}{l}
\varepsilon((0,0),id)+\varepsilon((1,0),id) = 0\\
\varepsilon((0,0),\perm{21})+\varepsilon((1,0),\perm{21}) = 0\\
\varepsilon((0,1),id)+\varepsilon((1,1),id) = 0\\
\varepsilon((0,1),\perm{21})+\varepsilon((1,1),\perm{21}) = 0
\end{array}
\right.$$
Solving these ensures that $\theta_{\varepsilon,\gamma}$ is a copula density.
To guarantee that, for any $g$ symmetric, $g(V)\eqd g(U)$, $\varepsilon$ must satisfy \eqref{epscondG}:
$$\left\{
\begin{array}{l}
\varepsilon((0,0),id)+\varepsilon((0,0),\perm{21}) = 0\\
\varepsilon((1,0),id)+\varepsilon((0,1),\perm{21}) = 0\\
\varepsilon((0,1),id)+\varepsilon((1,0),\perm{21}) = 0\\
\varepsilon((1,1),id)+\varepsilon((1,1),\perm{21}) = 0
\end{array}
\right.$$
We deduce that $\varepsilon$ must take the form
$$\varepsilon(\alpha,\beta) = \lambda(-1)^{|\alpha|}\sgn(\beta),$$
where $|\alpha|=\alpha_1+\alpha_2$ and $\sgn(\beta)$ is the signature of the permutation $\beta$, $\sgn(id)=1$ and $\sgn(\perm{21})=-1$.

In other words, \eqref{eqOctalCopDef} defines, up to a multiplicative factor, the only copula such that, for any $g\in\mathcal{G}_2$, $g(V)\eqd g(U)$.
\end{example}

Observe that if $\theta_{\varepsilon,\gamma}$ satisfies the conditions of the theorem, then so does $\theta_{\varepsilon/\lambda,\lambda\gamma}$, for any $\lambda>0$, and the same goes for $\theta_{-\varepsilon,\gamma}$ as long as
$$\gamma(u)\leq -\max_{\varepsilon(\alpha,\beta)<0}\frac{c(\tau_\alpha(\sigma_{\beta^{-1}}(u)))}{\varepsilon(\alpha,\beta)}.$$

While in general, \eqref{epscondM} and \eqref{epscondC} are sufficient to guarantee that marginal distributions match and that $\theta_{\varepsilon,\gamma}$ is a copula density, respectively, the next result shows that they may not be necessary, except if we add the assumption that $\gamma$ is bounded from above and away from 0.

\begin{theorem}\label{thmBoundCopula}
Further to the setting of Theorem \ref{thmGenCopula}, we assume that $\gamma$ is bounded from above and away from 0; i.e. we assume that
\begin{equation}\label{boundgamma}
\inf(\gamma)=\inf_{u\in\Delta(0,id)}\gamma(u)>0\mbox{ and }\sup(\gamma)=\sup_{u\in\Delta(0,id)}\gamma(u)<+\infty.
\end{equation}
\begin{enumerate}
\item[(M)] Fix $k\in[d]$. $(V_1,\ldots,V_{k-1},V_{k+1},\ldots,V_d) \eqd (U_1,\ldots,U_{k-1},U_{k+1},\ldots,U_d)$ if and only if
\begin{equation}\label{epscondMiff}
\forall j\in[d],\forall a\in\{0,1\}^{d-1},\forall b\in S_{d-1},\ \sum_{\aleph=0,1}\varepsilon(\omega_k(a,\aleph),\chi_j(b,k))=0.
\end{equation}
\item[(C)] Suppose $c$ is a copula density. $\theta_{\varepsilon,\gamma}$ is a copula density if and only if
\begin{equation}\label{epscondCiff}
\forall j,k\in[d],\ \sum_{a,b}\varepsilon(\omega_k(a,0),\chi_j(b,k))=\sum_{a,b}\varepsilon(\omega_k(a,1),\chi_j(b,k))=0.
\end{equation}
\end{enumerate}
\end{theorem}
\begin{proof}
The sufficiency of both statements was shown in Theorem \ref{thmGenCopula}. We limit the proofs to necessity.
\begin{enumerate}
\item[(M)] We know that $(V_1,\ldots,V_{k-1},V_{k+1},\ldots,V_d) \eqd (U_1,\ldots,U_{k-1},U_{k+1},\ldots,U_d)$ if and only if, $\forall a\in\{0,1\}^{d-1}$, $\forall b\in S_{d-1}$,
$$\forall q\in(0,1/2)^{d-1}\mbox{ such that }q_1<\ldots<q_{d-1},\ \sum_{j=1}^d\lambda_j\int_{[q_{j-1},q_j]}\gamma(\omega_j(q,r))dr=0,$$
where $\displaystyle\lambda_j = \sum_{\aleph=0,1}\varepsilon(\omega_k(a,\aleph),\chi_j(b,k))$, $q_0=0$ and $q_d=1/2$. We reason by contradiction and assume that $\{j\in[d]:\lambda_j\neq0\}\neq\emptyset$. Then
\begin{eqnarray*}
\lefteqn{\sum_{j=1}^d\lambda_j\int_{[q_{j-1},q_j]}\gamma(\omega_j(q,r))dr}\\
& \geq & \sum_{j:\lambda_j<0}\lambda_j\sup(\gamma)(q_j-q_{j-1}) + \sum_{j:\lambda_j>0}\lambda_j\inf(\gamma)(q_j-q_{j-1})
\end{eqnarray*}
and shrinking $q_j-q_{j-1}$ whenever $\lambda_j<0$ (and therefore expanding it whenever $\lambda_j>0$) shows that the right hand side can be made strictly positive, thus contradicting the assumption that the left hand side is nil. Of course, if there is no $j$ such that $\lambda_j>0$, then the inequality can be reversed and the left hand side shown to be strictly negative, leading to a contradiction.
\end{enumerate}
\item[(C)] We know that $\theta_{\varepsilon,\gamma}$ is a copula density if and only if $\forall k\in[d]$,
\begin{equation}\label{Ccond}
\forall r\in(0,1),\ \sum_{j=1}^d\lambda_j\int_{\Delta_j(0,id,r)}\gamma(\omega_j(v,\bar r))dv=0,
\end{equation}
where $\displaystyle\lambda_j=\sum_{a,b}\varepsilon(\alpha,\beta)$ and, $\alpha$ and $\beta$ are as in Lemma \ref{lemmarginal}. We also note that
$$\leb(\Delta_j(0,id,r)) = \frac{\bar r^{j-1}(1-2\bar r)^{d-j}}{2^{d-j}(j-1)!(d-j)!} \propto \bar r^{j-1}(1-2\bar r)^{d-j}.$$
Again we reason by contradiction. First we assume that $\lambda_1\neq0$ and more specifically (wlog) $\lambda_1>0$. Then
\begin{eqnarray*}
\lefteqn{\sum_{j=1}^d\lambda_j\int_{\Delta_j(0,id,r)}\gamma(\omega_j(v,\bar r))dv}\\
& \geq & \lambda_1\inf(\gamma)\leb(\Delta_1(0,id,r))+\sum_{j=2}^d\lambda_j\int_{\Delta_j(0,id,r)}\gamma(\omega_j(v,\bar r))dv.
\end{eqnarray*}
Since $\lim_{r\to0}\leb(\Delta_j(0,id,r))=0$, for $j\in\{2,\ldots,d\}$, and $\gamma$ is bounded, by making $r$ approach 0, the second term in the right hand side can be made as small as we want, while the first term is strictly positive. It follows that the left hand side can be made strictly positive thus contradicting the fact that it must be nil for all $r$. We deduce that $\lambda_1$ must be nil and \eqref{Ccond} becomes
$$\forall r\in(0,1),\ \sum_{j=2}^d\lambda_j\frac1{r}\int_{\Delta_j(0,id,r)}\gamma(\omega_j(v,\bar r))dv=0.$$
Again, we assume (wlog) that $\lambda_2>0$. Then
\begin{eqnarray*}
\lefteqn{\sum_{j=2}^d\lambda_j\frac1{r}\int_{\Delta_j(0,id,r)}\gamma(\omega_j(v,\bar r))dv}\\
& \geq & \lambda_2\inf(\gamma)\frac{\leb(\Delta_2(0,id,r))}{r}+\sum_{j=3}^d\lambda_j\frac1{r}\int_{\Delta_j(0,id,r)}\gamma(\omega_j(v,\bar r))dv
\end{eqnarray*}
and the second term in the right hand side can be made as small as we want, while the first term is strictly positive. It follows that the left hand side can be made strictly positive thus showing that $\lambda_2$ must be nil.
We continue this way, adjusting \eqref{Ccond} by increasing powers of $r$, to prove that $\lambda_1=\ldots=\lambda_{d-1}=0$ and finally that $\lambda_d=0$ since
$$\forall r\in(0,1),\ \lambda_d\int_{\Delta_d(0,id,r)}\gamma(\omega_d(v,\bar r))dv=0.$$
\end{proof}

\begin{corollary}
Suppose \eqref{boundgamma} holds and $\varepsilon$ takes the form
$$\varepsilon(\alpha,\beta) = \zeta(\alpha)\psi(\beta).$$
\begin{enumerate}
\item[(D$'$)] $\theta_{\varepsilon,\gamma}$ is a density if and only if
$$\text{either }\sum_\alpha\zeta(\alpha)=0\text{ or }\sum_\beta\psi(\beta)=0.$$
\item[(G$'$)] $g(V)\eqd g(U)$, for any $g$ symmetric, if and only if
$$\sum_{\beta}\psi(\beta)=0.$$
\item[(C$'$)] Suppose $c$ is a copula density. $\theta_{\varepsilon,\gamma}$ is a copula density if and only if
$$\forall k\in[d],\ \sum_a\zeta(\omega_k(a,0))=\sum_a\zeta(\omega_k(a,1))=0.$$
\end{enumerate}
\end{corollary}

\begin{proposition}
A necessary and sufficient condition for
$$\forall j,k\in[d],\forall a\in\{0,1\}^{d-1},\forall b\in S_{d-1},\ \sum_{\aleph=0,1}\varepsilon(\omega_k(a,\aleph),\chi_j(b,k))=0$$
is that
$$\forall\alpha\in\{0,1\}^d,\forall\beta\in S_d\ \varepsilon(\alpha,\beta)=(-1)^{|\alpha|}\varepsilon(0,\beta),$$
where $|\alpha|=\alpha_1+\ldots+\alpha_d$.
\end{proposition}
\begin{proof}
Sufficiency is immediate. We prove necessity by induction on $d$. First we observe that for any $\beta\in S_d$, for any $k\in[d]$, there exists $j\in[d]$ and $b\in S_{d-1}$ such that $\beta=\chi_j(b,k)$. Indeed, letting $j=\beta^{-1}(k)$ and
$$b(i)=\left\{
\begin{array}{ll}
\beta(i) & \text{if } i\leq j-1\mbox{ and }\beta(i)\leq k-1\\
\beta(i)-1 & \text{if } i\leq j-1\mbox{ and }\beta(i)\geq k+1\\
\beta(i+1) & \text{if } i\geq j\mbox{ and }\beta(i+1)\leq k-1\\
\beta(i-1)+1 & \text{if } i\geq j\mbox{ and }\beta(i+1)\geq k+1
\end{array}
\right.$$
we obtain the required identity. We shall therefore prove that, for any fixed $\beta\in S_d$, the condition
$$\forall k\in[d],\forall a\in\{0,1\}^{d-1},\ \varepsilon(\omega_k(a,0),\beta)+\varepsilon(\omega_k(a,1),\beta)=0$$
implies the desired statement. As $\beta$ is fixed throughout, we write $\zeta(\alpha)$ for $\varepsilon(\alpha,\beta)$.

The case $d=2$ can easily be checked. Suppose the necessity true for $d-1$. Then setting the first component in $\alpha$ and $a$ to 0 reduces the dimensionality of the problem by 1 and leads to
$$\forall\alpha\in\{0\}\times\{0,1\}^{d-1},\ \zeta(\alpha)=(-1)^{|\alpha|}\zeta(0).$$
Similarly, setting the first component in $\alpha$ and $a$ to 1 again reduces the dimensionality of the problem by 1 and leads to
$$\forall\alpha\in\{1\}\times\{0,1\}^{d-1},\ \zeta(\alpha)=(-1)^{|\alpha|-1}\zeta(\omega_1(0,1)).$$
Now taking $k=1$ and $a=0$ leads to $\zeta(\omega_1(0,1))=-\zeta(0)$ and concludes the proof.
\end{proof}

\begin{corollary}
Suppose $\displaystyle\varepsilon(\alpha,\beta) = (-1)^{|\alpha|}\psi(\beta)$.
If $\sum_\beta\psi(\beta)=0$ then all conditions of Theorem \ref{thmGenCopula} are satisfied; that is, for any
$\gamma$ such that
$$\gamma(u)\leq \min_{(-1)^{|\alpha|}\psi(\beta)>0}\frac{c(\tau_\alpha(\sigma_{\beta^{-1}}(u)))}{|\psi(\beta)|},$$
$\theta_{\varepsilon,\gamma}$ is a copula density for which the $(d-1)$-dimensional marginals coincide with those of $c$ and, for any $g$ symmetric, $g(V)\eqd g(U)$, where $U$ and $V$ have densities $c$ and $\theta_{\varepsilon,\gamma}$, respectively.

In particular, this is true for $\displaystyle\varepsilon(\alpha,\beta) = (-1)^{|\alpha|}\sgn(\beta)$,
where $\sgn(\beta)$ is the signature of the permutation $\beta$.
\end{corollary}

\begin{example}[The trivariate case]
When $d=3$, the problem offers sufficiently many degrees of freedom to permit multiple solutions. For example
$$\varepsilon(\alpha,\beta) = \zeta(\alpha)\sgn(\beta),$$
where
$\displaystyle\zeta(\alpha)=(-1)^{|\alpha|}1_{\alpha_3=0}$
provides a copula such that, for any $g\in\mathcal{G}_3$, $g(V)\eqd g(U)$.
However, in this case, the two-dimensional marginals do not match.
\end{example}

\subsection{The general multidimensional case}

We are now ready to deal with the case of $d$ identically distributed arbitrary random variables. We stress here that we do not assume that the random variables are independent.

\begin{proposition}
Suppose that $X_1,\ldots,X_d$ are identically distributed, that $X=(X_1,\ldots,X_d)$ has copula density $c$, marginal distribution function $\Phi$ and marginal density $\phi$, so that its density is
$$f(x) = c\big(\Phi(x_1),\ldots,\Phi(x_d)\big)\prod_{k=1}^d\phi(x_k).$$
For any $(\varepsilon,\gamma)$ satisfying conditions \eqref{epscondG} and \eqref{epscondC} of Theorem \ref{thmGenCopula},
$$\varphi(x) =\theta_{\varepsilon,\gamma}\big(\Phi(x_1),\ldots,\Phi(x_d)\big)\prod_{k=1}^d\phi(x_k)$$
generates a random variable $Y=(Y_1,\ldots,Y_d)$ that satisfies the requirements that, for any $k\in[d]$, $Y_k\eqd X_k$, and, for any $g\in\mathcal{G}_d$, $g(Y)\eqd g(X)$.
\end{proposition}
\begin{proof}
Let $U_k=\Phi(X_k)$, $g\in\mathcal{G}(\mathbb{R}^d)$ and
$$\mathfrak{g}(u)=g(\Phi^{-1}(u_1),\ldots,\Phi^{-1}(u_d)).$$
Then $U$ has density $c$ and $\mathfrak{g}\in\mathcal{G}([0,1]^d)$.
Letting $V$ be a random variable with density $\theta_{\varepsilon,\gamma}$, $Y_k=\Phi^{-1}(V_k)$, we get that
\begin{eqnarray*}
g(Y) & = & g(\Phi^{-1}(V_1),\ldots,\Phi^{-1}(V_d))\ =\ \mathfrak{g}(V)\\
& \eqd & \mathfrak{g}(U) = g(\Phi^{-1}(U_1),\ldots,\Phi^{-1}(U_d)) = g(X).
\end{eqnarray*}
\end{proof}

\begin{example}
Let $\Phi$ be the distribution function and $\phi$ be the density of the standard normal distribution. Then for any
$\gamma$ such that
$$\gamma(u)\leq \min_{(-1)^{|\alpha|}\sgn(\beta)>0}c(\tau_\alpha(\sigma_{\beta^{-1}}(u))),$$
$$\varphi(x) =\theta_\gamma\big(\Phi(x_1),\ldots,\Phi(x_d)\big)\prod_{k=1}^d\phi(x_k),$$
where
$$\theta_\gamma(u) = 1-\sum_{\alpha,\beta}(-1)^{|\alpha|}\sgn(\beta)\gamma(\sigma_\beta(\tau_\alpha(u)))1_{\Delta(\alpha,\beta)}(u),$$
is the density of a $d$-dimensional random variable $Y$ for which all $(d-1)$-dimensional marginals are independent and identically distributed standard normal random variables, $Y_1+\ldots+Y_d$ is normal with mean 0 and variance $d$, and $Y$ is non-Gaussian.
\end{example}

\subsection{The case of non-identically distributed random variables}\label{nonidentical}
Can the above construction extend to the case of non-identically distributed (and non-independent) random variables? To answer this question, we return to the two-dimensional case. Let $s_1$, $s_2$ and $s_{12}$ be the reflections
$$s_1(u_1,u_2)=(1-u_1,u_2),\ s_2(u_1,u_2)=(u_1,1-u_2)\mbox{ and }s_{12}(u_1,u_2)=(u_2,u_1).$$
These three involutions are such that $s_1s_2=s_2s_1$, $s_1s_{12}=s_{12}s_2$ and $s_2s_{12}=s_{12}s_1$. It follows that they generate a finite non-Abelian group $\mathcal{R}=\{id,s_1,s_2,s_{12},s_1s_2,s_1s_{12},s_2s_{12},s_1s_2s_{12}\}$, the dihedral group of order 8.

Each element of $\mathcal{R}$ corresponds to one of the eight regions $\Delta(\alpha,\beta)$, and $\varepsilon(\alpha,\beta)$ of \eqref{eqOctalCopDef} is simply $(-1)^{|s|}$, where $|s|$ is the word length of $s$, that is the number of generators in the decomposition of $s$ (modulo 2).

In the case of non-identically distributed random variables, say with distribution functions $\Phi_1$ and $\Phi_2$, for the construction to hold for symmetric functions, and in particular for the sum, the generator $s_{12}$ needs to be changed to
$$s_{12}(u_1,u_2) = (\Phi_1(\Phi_2^{-1}(u_2)),\Phi_2(\Phi_1^{-1}(u_1)).$$

One would then attempt a construction of the type
$$\theta(u) = c(u)-\sum_{s\in\mathcal{R}}\varepsilon(s)\gamma(s(u))1_\Delta(s(u)),$$
for an appropriate reference region $\Delta$, one for which $\{s(\Delta),s\in\mathcal{R}\}$ forms a measurable partition of $[0,1]^2$.
The next example illustrates the difficulty we face in general.

\begin{example}
Suppose $\Phi_2(x)=\Phi_1(x)^2$ ($X_2$ has the distribution of the maximum of two independent copies of $X_1$) so that $s_{12}(u_1,u_2)=(\sqrt{u_2},u_1^2)$. Then
$$s_1s_{12}s_2s_{12}(u_1,u_2) = \left(1-\sqrt{1-u_1^2},u_2\right)$$
and $(s_1s_{12}s_2s_{12})^n$, obtained by iterating the map $1-\sqrt{1-r^2}$, yields the infinite sequence
$$1-\sqrt{1-u_1^2},1-\sqrt{1-\left(1-\sqrt{1-u_1^2}\right)^2},1-\sqrt{1-\left(1-\sqrt{1-\left(1-\sqrt{1-u_1^2}\right)^2}\right)^2},\ldots.$$
In this case $\mathcal{R}$ is infinite and it is not at all clear how to find $\Delta$ (assuming it exists).
\end{example}

In the above example the identity $s_1s_{12}=s_{12}s_2$ fails resulting in $\mathcal{R}$ being infinite. This identity translates in the language of the previous subsection to $\displaystyle\sigma_\beta(\tau_\alpha(u)) = \tau_{\sigma_\beta(\alpha)}(\sigma_\beta(u))$ which was crucial in our construction. It enabled us to use the symmetry of $g$ in the proof of (G) (see Theorem \ref{thmGenCopula}) to establish that
\begin{eqnarray*}
\int g(u)\gamma(s_{12}s_1(u))1_\Delta(s_{12}s_1(u))du & = & \int_\Delta g(s_1s_{12}(u))\gamma(u)du\ =\ \int_\Delta g(s_{12}s_2(u))\gamma(u)du\\
& = & \int_\Delta g(s_2(u))\gamma(u)du\ =\ \int g(u)\gamma(s_2(u))1_\Delta(s_2(u))du,
\end{eqnarray*}
and consequently obtain the requirement that $\varepsilon(s_{12}s_1)+\varepsilon(s_2)=\varepsilon((1,0),\perm{21})+\varepsilon((0,1),id)=0$.

\begin{figure}[!htb]
\begin{center}
\includegraphics[height=5cm]{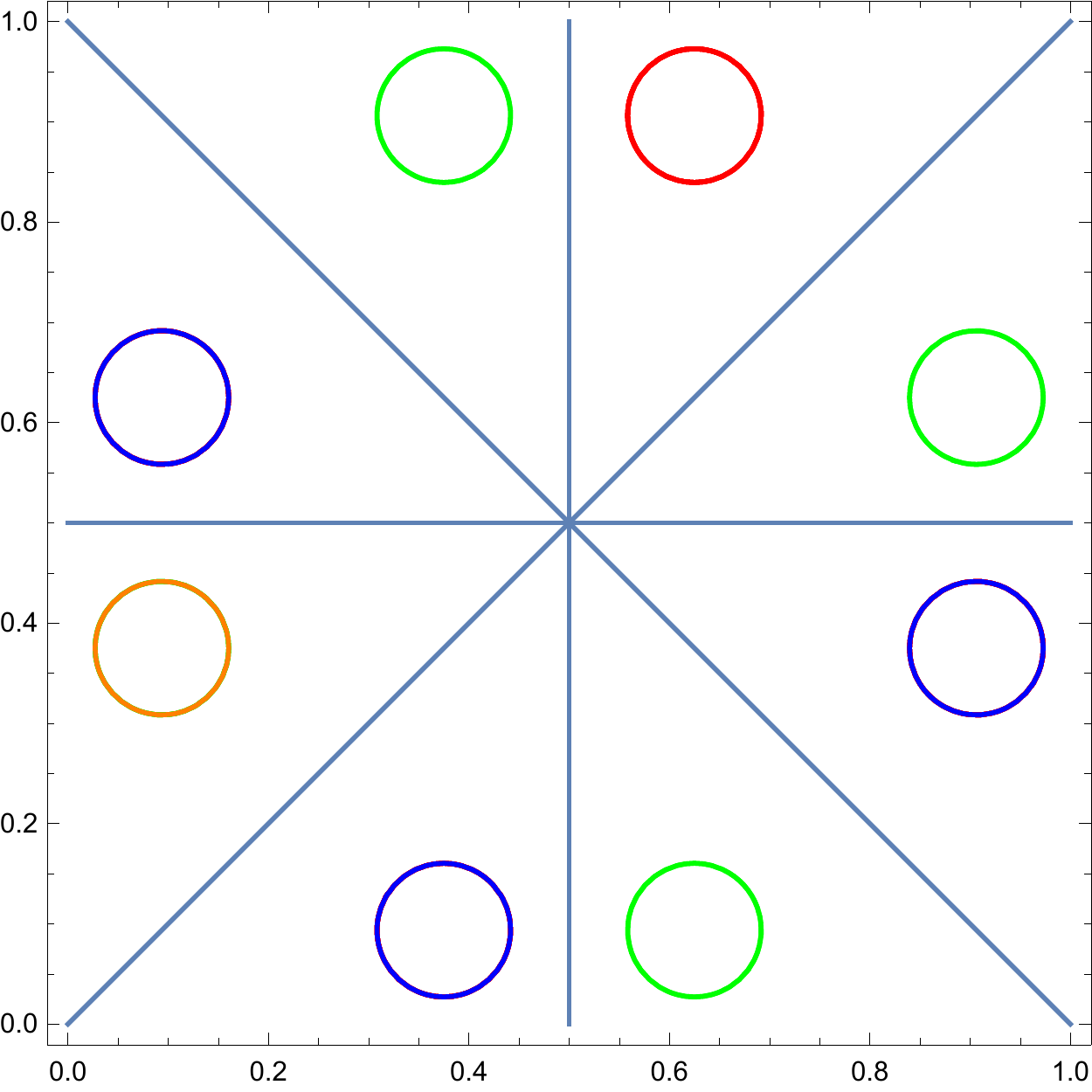}\includegraphics[height=5cm]{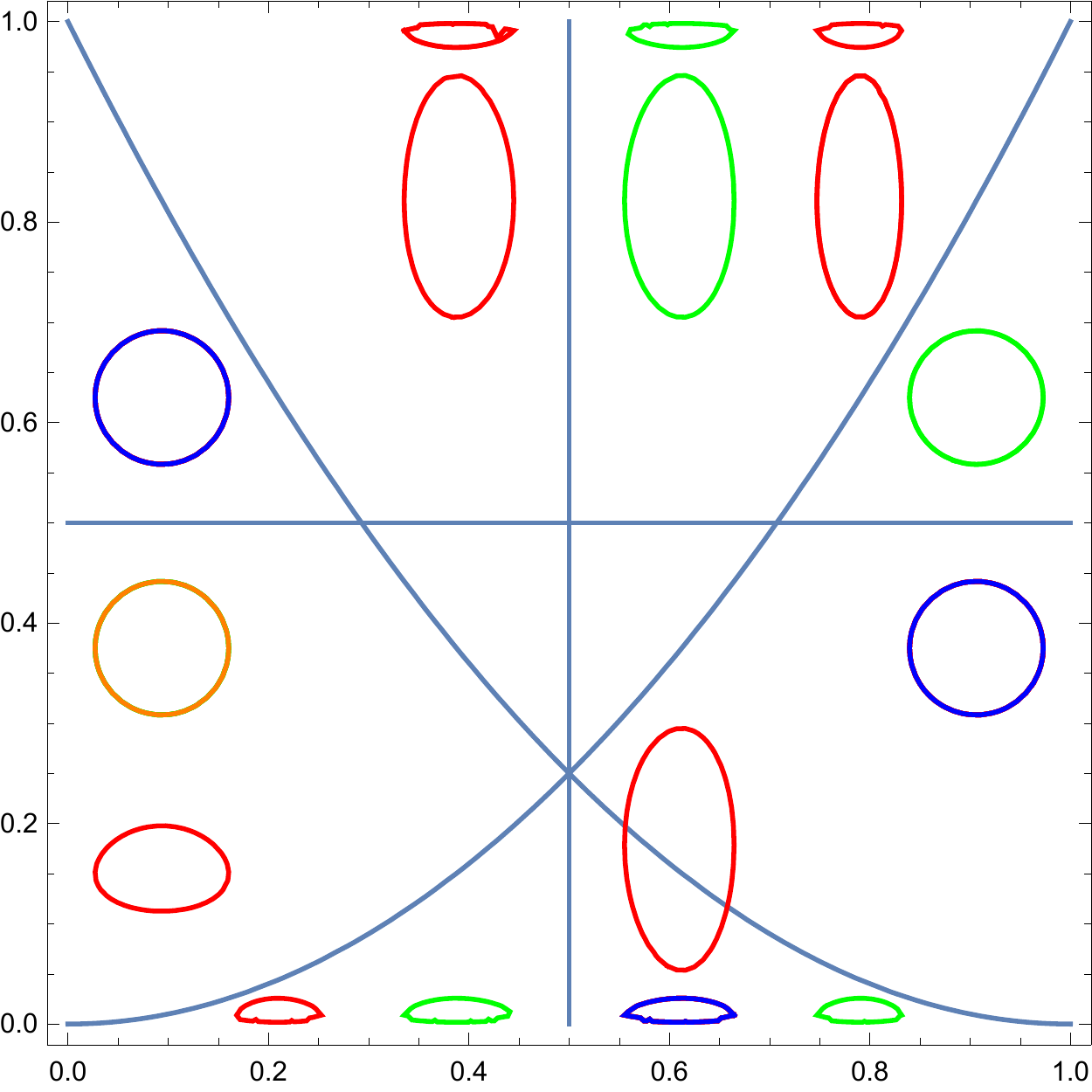}\includegraphics[height=5cm]{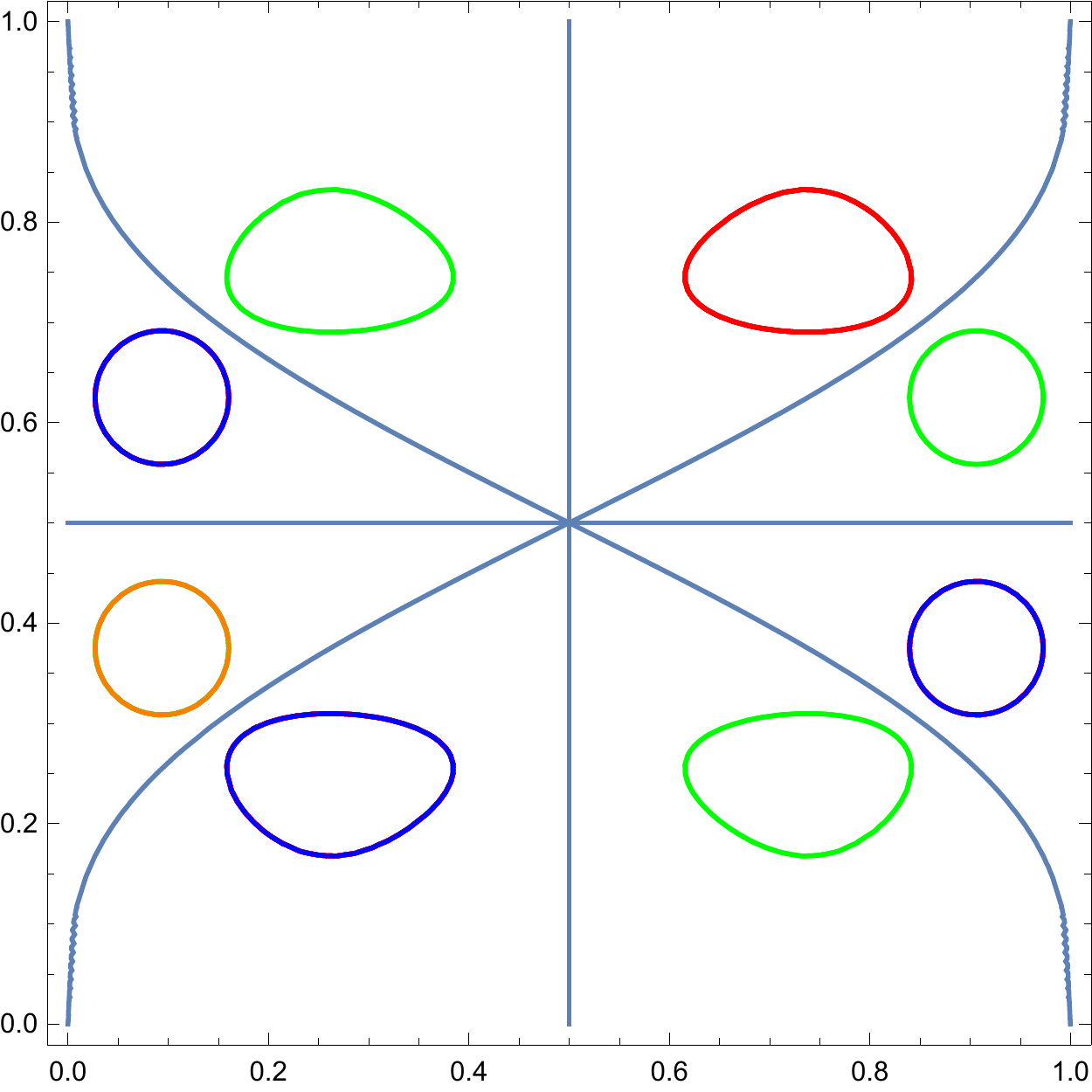}
\end{center}
\caption{The images of a circle (gold) by all elements of $\mathcal{R}$ of length 1 (purple), 2 (green) and 3 (red), in the identically distributed case (left), the general non-identically distributed case (middle) and the similarly distributed case (right)}
\label{NonIID}
\end{figure}

In order to retain the identity $s_1s_{12}=s_{12}s_2$ we introduce the following notion -- see Figure \ref{NonIID} below for the motivation behind it.

\begin{definition*}
Two random variables are said to be similarly distributed if their distribution functions $\Phi_1$ and $\Phi_2$, assumed to be continuous and strictly increasing (on some interval), satisfy the identity
$$\Phi_1^{-1}(1-\Phi_1(x)) = \Phi_2^{-1}(1-\Phi_2(x)).$$
\end{definition*}

\begin{proposition}
\begin{enumerate}
\item Two identically distributed random variables are necessarily similarly distributed.
\item Two similarly distributed random variables must have the same median.
\item Two symmetrical distributions around the same median are similarly distributed.
In particular two normal random variables with the same mean (but different variances) are similarly distributed.
\end{enumerate}
\end{proposition}
\begin{proof}
(1) is trivial.

Let $m$ be the median for $\Phi_1$. Then
$$m = \Phi_1^{-1}(1-\Phi_1(m)) = \Phi_2^{-1}(1-\Phi_2(m))\mbox{ and }\Phi_2(m) = 1-\Phi_2(m),$$
showing that $m$ is also the median for $\Phi_2$.

Finally, if $\Phi_1$ and $\Phi_2$ are symmetric about the same median $m$, then
$$\Phi_1^{-1}(1-\Phi_1(x)) = m-x = \Phi_2^{-1}(1-\Phi_2(x)).$$
\end{proof}

For two similarly distributed random variables with distribution functions $\Phi_1$ and $\Phi_2$, we define
\begin{equation}\label{generators}
\sigma_1(x_1,x_2) = (\Phi_1^{-1}(1-\Phi_1(x_1)),x_2),\ \sigma_2(x_1,x_2) = (x_1,\Phi_2^{-1}(1-\Phi_2(x_2)))\mbox{ and }\sigma_{12}(x_1,x_2) = (x_2,x_1).
\end{equation}

\begin{proposition}
The three involutions $\sigma_1$, $\sigma_2$ and $\sigma_{12}$ are such that $\sigma_1\sigma_2=\sigma_2\sigma_1$, $\sigma_1\sigma_{12}=\sigma_{12}\sigma_2$ and $\sigma_2\sigma_{12}=\sigma_{12}\sigma_1$. As such, they generate a finite group $\mathcal{R}=\{id,\sigma_1,\sigma_2,\sigma_{12},\sigma_1\sigma_2,\sigma_1\sigma_{12},\sigma_2\sigma_{12},\sigma_1\sigma_2\sigma_{12}\}$.
\end{proposition}
\begin{proof}
The facts that $\sigma_1$, $\sigma_2$ and $\sigma_{12}$ are involutions and that $\sigma_1\sigma_2=\sigma_2\sigma_1$ are easily checked. We prove that $\sigma_1\sigma_{12}=\sigma_{12}\sigma_2$; $\sigma_2\sigma_{12}=\sigma_{12}\sigma_1$ follows by taking inverses:
\begin{eqnarray*}
\sigma_1\sigma_{12}(x_1,x_2) & = & \sigma_1(x_2,x_1)\ =\ (\Phi_1^{-1}(1-\Phi_1(x_2)),x_1)\ =\ (\Phi_2^{-1}(1-\Phi_2(x_2)),x_1)\\
& = & \sigma_{12}(x_1,\Phi_2^{-1}(1-\Phi_2(x_2)))\ =\ \sigma_{12}\sigma_2(x_1,x_2).
\end{eqnarray*}
\end{proof}

While it is possible to approach this situation via copulas, other than in the identically distributed case, the resulting $\theta$ turns out to depend on $\Phi_1$ and $\Phi_2$ making it not universal and therefore less desirable. Instead, we apply the symmetry-balancing approach directly to the density.
\begin{theorem}
Let $f$ be the joint density of two similarly distributed random variables, $X_1$ and $X_2$, $m$ be the common median, $\Delta=\{x\in\mathbb{R}^2:x_1<x_2<m\}$ and $(\varepsilon,\gamma)$ be such that
\begin{equation}\label{varphi}
\varphi(x_1,x_2) = f(x_1,x_2)-\sum_{\sigma\in\mathcal{R}}(-1)^{|\sigma|}\gamma(\sigma(x_1,x_2))|J_\sigma(x_1,x_2)|1_\Delta(\sigma(x_1,x_2))
\end{equation}
is non-negative, where $J_\sigma$ denotes the Jacobian determinant of $\sigma$ and $|\sigma|$ is the word length of $\sigma$, that is the number of generators in the decomposition of $\sigma$ (modulo 2).

Then $\varphi$ generates $(Y_1,Y_2)$ such that $Y_1\eqd X_1$, $Y_2\eqd X_2$ and, for any $g\in\mathcal{G}_2$, $g(Y_1,Y_2)\eqd g(X_1,X_2)$; in particular $Y_1+Y_2\eqd X_1+X_2$.
\end{theorem}
\begin{proof}
Let $\Psi(x)=\Phi_1^{-1}(1-\Phi_1(x)) = \Phi_2^{-1}(1-\Phi_2(x))$ and $\psi(x)=\Psi'(x)$. We start by checking that $\int\varphi(x_1,x_2)dx_1=\int f(x_1,x_2)dx_1$. The integration in $x_2$ is performed in an identical manner. Suppose $x_2<m$. Then $\Psi(x_2)>m$ and $\sigma_2(x_1,x_2)=(x_1,\Psi(x_2))\not\in\Delta$, $\sigma_{12}\sigma_2(x_1,x_2))=(\Psi(x_2),x_1)\not\in\Delta$, $\sigma_1\sigma_2\sigma_{12}(x_1,x_2))=(\Psi(x_2),\Psi(x_1))\not\in\Delta$ and $\sigma_1\sigma_2(\Psi(x_1),\Psi(x_2)))=(\Psi(x_2),x_1)\not\in\Delta$. It follows that the sum in \eqref{varphi} only contains four non-zero expressions, those for $\sigma\in\{id,\sigma_1,\sigma_{12},\sigma_{12}\sigma_1\}$.
Now,
\begin{eqnarray*}
\lefteqn{\sum_{\sigma\in\{id,\sigma_1\}}(-1)^{|\sigma|}\int\gamma(\sigma(x_1,x_2))|J_\sigma(x_1,x_2)|1_\Delta(\sigma(x_1,x_2))dx_1}\\
& = & \int\gamma(x_1,x_2)1_\Delta(x_1,x_2)dx_1-\int\gamma(\sigma_1(x_1,x_2))|J_{\sigma_1}(x_1,x_2)|1_\Delta(\sigma_1(x_1,x_2))dx_1\\
& = & \int_{(-\infty,x_2)}\gamma(x_1,x_2)dx_1-\int_{(\Psi(x_2),+\infty)}\gamma(\Psi(x_1),x_2)|\psi(x_1)|dx_1\\
& = & \int_{(-\infty,x_2)}\gamma(x_1,x_2)dx_1-\int_{(-\infty,x_2)}\gamma(z_1,x_2)dz_1\ =\ 0
\end{eqnarray*}
and
\begin{eqnarray*}
\lefteqn{\sum_{\sigma\in\{\sigma_{12},\sigma_{12}\sigma_1\}}(-1)^{|\sigma|}\int\gamma(\sigma(x_1,x_2))|J_\sigma(x_1,x_2)|1_\Delta(\sigma(x_1,x_2))dx_1}\\
& = & -\int\gamma(x_2,x_1)1_\Delta(x_2,x_1)dx_1+\int\gamma(\sigma_{12}\sigma_1(x_1,x_2))|J_{\sigma_{12}\sigma_1}(x_1,x_2)|1_\Delta(\sigma_{12}\sigma_1(x_1,x_2))dx_1\\
& = & -\int_{(x_2,m)}\gamma(x_2,x_1)dx_1+\int_{(m,\Psi(x_2))}\gamma(x_2,\Psi(x_1))|\psi(x_1)|dx_1\\
& = & -\int_{(x_2,m)}\gamma(x_2,x_1)dx_1+\int_{(x_2,m)}\gamma(x_2,z_1)dz_1\ =\ 0.
\end{eqnarray*}
Similar calculations can be performed when $x_1$ and $x_2$ are swaped. In all, $\varphi$ and $f$ have the sane marginal distributions.

Furthermore, if $g$ is bounded and symmetric, then
\begin{eqnarray*}
\lefteqn{\sum_{\sigma\in\{id,\sigma_{12}\}}(-1)^{|\sigma|}\int g(x)\gamma(\sigma(x))|J_\sigma(x)|1_\Delta(\sigma(x))dx}\\
& = & \int_\Delta g(x)\gamma(x)dx-\int_\Delta g(\sigma_{12}(x))\gamma(x)dx\ =\ \int_\Delta g(x)\gamma(x)dx-\int_\Delta g(x)\gamma(x)dx\ =\ 0,
\end{eqnarray*}
\begin{eqnarray*}
\lefteqn{\sum_{\sigma\in\{\sigma_2,\sigma_2\sigma_{12}\}}(-1)^{|\sigma|}\int g(x)\gamma(\sigma(x))|J_\sigma(x)|1_\Delta(\sigma(x))dx}\\
& = & -\int g(x)\gamma(\sigma_2(x))|\psi(x_2)|1_\Delta(\sigma_2(x))dx+\int g(x)\gamma(\sigma_2\sigma_{12}(x))|\psi(x_1)|1_\Delta(\sigma_2\sigma_{12}(x))dx\\
& = & -\int_\Delta g(\sigma_2(z))\gamma(z)dz+\int g(\sigma_{12}\sigma_2(z))\gamma(z)dz\ =\ -\int_\Delta g(\sigma_2(z))\gamma(z)dz+\int g(\sigma_2(z))\gamma(z)dz\ =\ 0,
\end{eqnarray*}
and so on for the sums on $\{\sigma_1,\sigma_1\sigma_{12}\}$ and $\{\sigma_1\sigma_2,\sigma_1\sigma_2\sigma_{12}\}$, which concludes the proof.
\end{proof}

%Denote $D= [d]\times [d]-\{(i,i),i\in [d]$ as the off-diagonal elements of
%$[d]\times [d]$. Then the corresponding density is
%\begin{eqnarray}
%&&\prod_{j=1}^d\frac{1}{\sqrt{2\pi}\sigma_j}	e^{-\frac{1}{2}(x_j^2/\sigma_j^2)}
%\nonumber \\
%&&+ \kappa\prod_{j=1}^d\frac{1}{\sqrt{2\pi}\eta}	e^{-\frac{1}{2}(x_j^2/\eta^2)}
%\nonumber \\
%&&\times
%\sum_{A \subset D}(-2)^{|A|}\prod_{(i,j)\in A}H_2(x_i;\eta)H_2(x_j;\eta)\prod_{(k,l) \in D-A}\big ((H_3(x_k;\eta)+H_3(x_l;\eta)\big )
%\nonumber \\
%\label{inversion:0}
%\end{eqnarray}
%

\end{document}